\documentclass[11pt]{amsart}

\usepackage[USenglish]{babel}

\usepackage{latexsym,amssymb,amsmath,youngtab}\usepackage{multirow}
\usepackage{youngtab}
\usepackage{epsfig}
\usepackage{tikz,parskip}\usepackage[margin=2cm]{geometry}
\usepackage{tikz}
\usetikzlibrary{arrows,shapes.geometric,positioning,decorations.markings}

\usepackage{amsmath,amsthm,amssymb,amscd,MnSymbol}
\usepackage[mathscr]{eucal}
\usepackage{color}
\newcommand*{\textlabel}[2]{%
  \edef\@currentlabel{#1}
  \phantomsection
  #1\label{#2}
}

\newtheorem{theorem}{Theorem}[section]
\newtheorem{proposition}[theorem]{Proposition}

\newtheorem{proposition/definition}[theorem]{Proposition/Definition}

\theoremstyle{definition}
\newtheorem{definition}[theorem]{Definition}

\newtheorem{example}[theorem]{Example}

\theoremstyle{remark}
\newtheorem{remark}[theorem]{Remark}

\newcommand{\stack}[2]{\ensuremath{\genfrac{}{}{0pt}{}{#1}{#2}}} 

\newtheoremstyle{exercise}
  {3pt}
  {6pt}
  {}
  {}
  {\bfseries}
  {:}
  { }
   {}
\theoremstyle{exercise}
\newtheorem{exercise}[theorem]{Exercise}
\newtheoremstyle{exercises}
  {3pt}
  {6pt}
  {}
  {}
  {\bfseries}
  {:}
  {\newline}
   {}

\theoremstyle{exercise}
\newtheorem{exercises}[theorem]{Exercises}

\input epsf
\def\boxit#1{\vbox{\hrule height1pt\hbox{\vrule width1pt\kern3pt
  \vbox{\kern3pt#1\kern3pt}\kern3pt\vrule width1pt}\hrule height1pt}}

\def\cD{{\mathcal D}}

\def\av{av}

\def\tsup{\operatorname{sup}}\def\tinf{\operatorname{inf}}

\def\trank{\text{rank}}

\def\bee{\bold e}

\def\BC{\mathbb C}\def\BF{\mathbb F}
\def\BR{\mathbb R}
\def\BP{\mathbb P}
\def\pp#1{\mathbb P^{#1}}

\def\pp#1{{\mathbb P}^{#1}}
\def\tdim{{\rm dim}}

\def\hd{,...,}
\def\ww{\wedge}

\def\inv{{}^{-1}}

\def\tmult{{\rm mult}}
\def\cX{{\mathcal X}}\def\cY{{\mathcal Y}}
\def\cZ{{\mathcal Z}}\def\cM{{\mathcal M}}
\def\cH{{\mathcal H}}
\def\cE{{\mathcal E}}

\def\cS{{\mathcal S}}

\def\11{\mathbf 1}

\def\fu{{\mathfrak u}}

\def\l{\lambda}
\def\a{\alpha}

\def\o{\omega}

\def\b{\beta}

\def\s{\sigma}

\def\d{\delta}
\def\th{\theta}

\def\ot{{\mathord{ \otimes } }}
\def\op{{\mathord{\,\oplus }\,}}
\def\otc{{\mathord{\otimes\cdots\otimes} }}

\def\ra{{\mathord{\;\rightarrow\;}}}

\def\La#1{\Lambda^{#1}}

\def\op{\oplus}
\def\BF{\Bbb F}\def\BZ{\Bbb Z}

\def\ep{\epsilon}
\def\op{\oplus}

\def\s{\sigma}

\def\a{\alpha}
\def\b{\beta}

\def\l{\lambda}

\def\FS{\mathfrak  S}

\def\ol{\overline}

\def\BP{\mathbb  P}
\def\BC{\mathbb  C}

\def\pp#1{\mathbb  P^{#1}}

\def\BR{\mathbb  R}

\def\tcos{\cos}
\def\tsin{\sin}

\def\ep{\epsilon}

\def\hd{, \hdots ,}

\def\inv{{}^{-1}}

\def\La#1{\Lambda^{#1}}

\def\pp#1{\mathbb  P^{#1}}

\def\ra{\rightarrow}

\def\ttrace{\operatorname{trace}}
\def\tend{\operatorname{End}}

\def\tdim{\operatorname{dim}}

\def\tlim{\lim}
\def\tmod{\operatorname{mod}}

\def\tmin{\operatorname{min}}
\def\tmax{\operatorname{max}}

\def\trank{\operatorname{rank}}

\def\ww{\wedge}

\def\ctimes{\times \cdots\times}

\def\bU{{\bold{U}}}

\def\be{\begin{equation}}
\def\ene{\end{equation}}

\DeclareMathOperator{\tlog}{log}
\DeclareMathOperator\tspan{span}\DeclareMathOperator\tsupp{supp}
\DeclareMathOperator\capp{Cap}
\DeclareMathOperator\trate{Rate}\DeclareMathOperator\tspec{spec}
\DeclareMathOperator\tsingvals{singvals}

\def\p{{\bold P}}
\def\np{{\bold N\bold P}}

\def\G{\Gamma}

\def\tspan{{\rm span}}

\newcommand{\Id}{\operatorname{Id}}

\newcommand{\tconv}{\operatorname{conv}}

\def\bra#1{|{#1}\rangle}\def\ket#1{\langle {#1}|}
\def\braket#1#2{\langle {#1}|{#2}\rangle}  
\def\ketbra#1#2{ \bra {#1}\ket {#2}}

\def\epr{\bra{epr}}

\def\eprx{\frac 1{\sqrt 2}(\bra{00}+\bra{11})}

\def\bra#1{|{#1}\rangle}\def\ket#1{\langle {#1}|}
\def\braket#1#2{\langle {#1}|{#2}\rangle}  
\def\ketbra#1#2{ \bra {#1}\ket {#2}}
\def\bU{{\bold{U}}}

\def\tinf{{\rm inf}}

\def\trank{{\mathrm {rank}}}

\def\tmult{{\rm mult}}

\begin{document}

\title[Quantum computation and information for mathematicians]{A very brief introduction to quantum computing and quantum information theory for mathematicians}
\author{J.~M. Landsberg}
\thanks{Landsberg    supported by   NSF grant  DMS-1405348.}
\begin{abstract}This is a very brief introduction to quantum computing and quantum information theory, primarily
aimed at geometers. Beyond basic definitions and examples, I emphasize aspects of interest to geometers, especially  
connections with asymptotic representation theory.
Proofs of most statements can be found in standard references such as \cite{MR1907291,MR1796805} as well as \cite{Lquantumnotes}.
\end{abstract}
\date\today
\maketitle

\section{Overview}
The postulates of quantum mechanics are presented as a generalized probability theory in \S\ref{qmasprob}. In \S\ref{EP} I
describe basic entanglement phenomena of \lq\lq super dense coding\rq\rq , \lq\lq teleportation\rq\rq, and Bell's confirmation of the \lq\lq paradox\rq\rq\
proposed by Einstein-Podolsky-Rosen. In \S\ref{Qalgs}  I outline aspects of the basic quantum algorithms, emphasizing the geometry involved.
\S\ref{cinfo} is a detour into classical information theory, which is the basis of its quantum cousin briefly  discussed in \S\ref{quantumchan}.
Before that, in \S\ref{quantumreform}, I reformulate quantum theory in terms of density operators, which facilitates the description
of quantum information theory. Critical to quantum information theory  is {\it von Neumann entropy} and in \S\ref{vnplus} I elaborate
on some of its properties. A generalization of \lq\lq teleportation\rq\rq\  (LOCC) is discussed in \S\ref{ELOCC}. I briefly discuss
SLOCC  in \S\ref{slocc} because of its geometric appeal. Regarding practical computation, the exponential growth in 
size of $(\BC^2)^{\ot n}$ with $n$ that appears in quantum information theory leads
to the notion of \lq\lq feasible\rq\rq\ states discussed in \S\ref{tns}, which has interesting
algebraic geometry associated to it. I conclude with a discussion of representation-theoretic aspects of quantum
information theory, including a discussion of the quantum marginal problem in \S\ref{rtinqt}.
I do not discuss topological quantum computing, which utilizes the  representation theory of the braid group.

\section{Quantum computation as generalized probabilistic computation}\label{qmasprob}
\subsection{Classical and Probabilistic computing via linear algebra}
\label{problinalg} 
This section is inspired by \cite[Exercise 10.4]{MR2500087}.

Classical communication deals with {\it bits}, elements of $\{ 0,1\}$, which will be convenient
to think of as elements of 
  $\BF_2$, the field with two elements. Let $f_n: \BF_2^n\ra \BF_2$ be a sequence of   functions.
  Give $\BR^2$ basis $\{\bra 0,\bra 1\}$ (such notation is standard in quantum mechanics) and give  $(\BR^2)^{\ot m}=\BR^{2^m}$
basis $\{\bra I\mid I\in \{0,1\}^m\}$. In this way, we may identify $\BF_2^m$ with the set of basis vectors of
$\BR^{2^m}$. 
A computation   of $f_n$ (via an arithmetic or Boolean circuit)  may be phrased as a sequence of   linear maps 
on a vector space containing $\BR^{2^n}$, where each linear map comes from a pre-fixed set agreed upon in advance.
In anticipation of what will come in quantum computation,
the pre-fixed set of maps will be taken from maps having the following properties:

\begin{enumerate}
\item 
Each linear map must take probability distributions to probability distributions. This implies the 
matrices are {\it stochastic}: the entries are non-negative and   each column sums to $1$.
\item
Each linear map   only alters a small number of 
  entries. For simplicity    
assume it alters at most three entries, i.e., it acts on at most   $\BR^{2^3}$ and is the identity on all other
factors in the tensor product. 
\end{enumerate}

To facilitate comparison with quantum comptuation, we will restrict to  reversible 
classical computation. The complexity class of a sequence of functions in  classical reversible computation is the same
as in arbitrary classical computation. 

 For example, if we want to effect $(x,y)\mapsto x*y$,
 consider the map 
\be\label{tofmap}
\bra{ x,y,z}\mapsto \bra{x,y,z\op (x*y)}=\bra{ x,y,z\op (x\ww y)}
\ene
(where the second expression is for those preferring Boolean notation)
and act as the identity on all other basis vectors (sometimes
called {\it registers}). Here $z$ will represent \lq\lq workspace  bits\rq\rq:  $x,y$ will come from 
the input  and $z$ will always be set to $0$ in the input. 
In the basis $\bra{000}, \bra{001},\bra{010},\bra{100},\bra{011},\bra{101},\bra{110},\bra{111} $, 
of $\BR^8$, 
  the matrix  is 
\be\label{tofmat}
\begin{pmatrix} 
1 &0 &0 &0 &0 &0 &0&0\\
0 &1 &0 &0 &0 &0 &0&0\\
0 &0 &1 &0 &0 &0 &0&0\\
0 &0 &0 &1 &0 &0 &0&0\\
0 &0 &0 &0 &1 &0 &0&0\\
0 &0 &0 &0 &0 &1 &0&0\\
0 &0 &0 &0 &0 &0 &0&1\\
0 &0 &0 &0 &0 &0 &1&0
\end{pmatrix}.
\ene
This gate is sometimes called the {\it Toffoli gate} and the matrix the {\it Toffoli matrix}.

The swap (negation) gate $\neg$ may be effected 
by the matrix
\be\label{negationmat}
\s_x=\begin{pmatrix} 0 & 1\\ 1&0\end{pmatrix}.
\ene
The swap and Toffoli matrix can perform any computation that is accomplished via a
sequence of matrices drawn from some finite set of Boolean operations,  each acting on a fixed number of basis vectors   with
at worst a polynomial  in $n$
size increase  in the number of matrices   needed.
For those familiar with  Boolean circuits, any sequence of Boolean circuits (one for each $n$) may be replaced by a sequence with just Toffoli and negation
gates with at worst a polynomial (in $n$)  blow up in size.

A probability distribution
on $\{0,1\}^{m}$ may be encoded  as a vector in $\BR^{2^m}$:
  If the probability distribution
assigns probability $p_I$ to 
  $  I\in \{0,1\}^{m}$, assign  to the distribution   the vector
$v=\sum_I p_I \bra I\in \BR^{2^m}$.

The above matrices realize  classical computation. To add randomness to enable probabilistic computation, introduce
the matrix
$$
\begin{pmatrix}
\frac 12 & \frac 12\\ \frac 12 &\frac 12\end{pmatrix}
$$
which acts on a single $\BR^2$ corresponding to a fair coin flip.
Note that the coin flip matrix is not invertible, which will be one motivation for
quantum computation in \S\ref{quantumotivate}.
Work in $\BR^{2^{n+s+r}}$ where $r$ is the number of times one needs  to access
a random choice and $s$ is the number of matrices (arithmetic operations)
in addition to the coin tosses needed to compute $f$.

A probabilistic computation, viewed this way, starts with $\bra{x0^{r+s}}$,
where $x\in \BF_2^n$ is the input. One then  applies a sequence of admissible
stochastic   linear maps  to
it, and ends with a vector that encodes a probability distribution on   $\{ 0,1\}^{n+s+r}$.
One then restricts this to $\{0,1\}^{p(n)}$, that is, one takes the vector and throws away
all but the first $p(n)$ entries. This vector encodes a probability sub-distribution, i.e.,
all coefficients are non-negative and they sum to a number between zero and one. 
One then renormalizes (dividing each entry by the sum of the entries) to
obtain a vector encoding a probability distribution on $\{0,1\}^{p(n)}$ and
then outputs the answer according to this distribution.
Note that  even if  our calculation is feasible (i.e., polynomial
in size), to write out the original output vector that one truncates would be exponential
in cost.  A stronger variant of this phenomenon will occur  with quantum computing, where the result
will be  obtained with a polynomial size calculation, but one does not have access to the vector
created, even using an exponential amount of computation.

To further prepare for the analogy with quantum computation,
define a probabilistic bit (a  {\it pbit})  to be an element of  
$$\{ p_0 \bra 0 +p_1 \bra 1 \mid p_j\in [0,1]{\rm \ and\ } p_0+p_1=1\}
\subset \BR^2.
$$

Note that the set of  pbits  (possible states)  is a convex set, and the basis
vectors are the extremal points of this convex set. 

\subsection{A wish list}\label{quantumotivate}
  Here is a wish list for how one might want to improve upon the above set-up:
\begin{enumerate}
\item Allow more general kinds of linear maps to get more computing power, while keeping the maps easy to compute.
\item Have reversible computation: we saw that classical computatation can be made reversible, but the coin flip
was not. This property is motivated by physics, where many physical theories require time reversibility.
\item Again motivated by physics, one would like to have a continous evolution of the probability vector, more
precisely, one would like the probability vector to depend on a continuous parameter $t$ such that  if
$\bra{\psi_{t_1}}=X\bra{\psi_{t_0}}$, then there exist admissible matrices $Y,Z$
such that $\bra{\psi_{t_0+\frac 12t_1}}=Y\bra{\psi_{t_0}}$ and $\bra{\psi_{t_1}}=Z\bra{\psi_{t_0+\frac 12 t_1}}$
and $X=ZY$. In particular, one wants operators to have square roots.
\end{enumerate}

One way to make the coin flip reversible is, instead of making the probability distribution be
determined by the sum of the coefficients, one could take the sum of the squares. If one does
this, there is no harm in allowing the entries of the output vectors to become negative,
and   one could
use 
$$
H:=\frac 1{\sqrt{2}}\begin{pmatrix} 1&1\\ 1& -1\end{pmatrix}
$$
for the coin flip. The matrix $H$ is called the {\it Hadamard matrix} or {\it Hadamard gate}
in the quantum computing literature.  
If we make this change, we obtain our second wish, and moreover have many operations be \lq\lq continous\rq\rq, 
because the set of matrices preserving the norm-squared of a real-valued vector is the  {\it  orthogonal
group}  $O(n)=\{ A\in Mat_{n\times n} \mid AA^T=\Id\}$. So for example, any rotation has a square root.

However our third property will not be completely satisfied, as the matrix
$$
\begin{pmatrix} 1&0\\ 0&-1\end{pmatrix}
$$
which represents a reflection, does not have a square root in $O(2)$.

To have the third wish satisfied,   allow  vectors with {\it complex} entries.
From now on   let $i=\sqrt{-1}$. For a complex number $z=x+iy$  let 
$\ol{z}=x-iy$ denote its complex conjugate and $|z|^2=z\ol{z}$ the square of its norm.

So we go from pbits, $\{p\bra 0 +q\bra 1\mid  \ p,q\geq 0\ {\rm{and}} \ p+q=1\}$
to 
{\it qubits}, the set of which is
\be\label{qubitdef}
\{\a \bra 0 +\b \bra 1\mid  \a,\b\in \BC\ {\rm and}\ |\a|^2+|\b|^2=1 \}.
\ene

The set of qubits, considered in terms of real parameters, looks at first
like the $3$-sphere $S^3$ in $\BR^4\simeq\BC^2$. However, the probability distributions
induced by   $\bra \psi$ and $e^{i\th}\bra\psi$ are the same so it is
really $S^3/S^1$ (the Hopf fibration), i.e., the two-sphere $S^2$.
Geometrically, it would be more natural (especially since we have already seen
the need to re-normalize in probabilistic computation) to work with projective space $\BC\pp 1\simeq S^2$
as our space of qubits, instead of a subset of $\BC^2$. So the set of qubits is better seen as \eqref{qubitdef} modulo
the equivalence $\bra\psi\sim e^{i\th}\bra\psi$.

For $v=(v_1\hd v_n)\in \BC^n$, write 
$|v|^2=|v_1|^2+\cdots +|v_n|^2$. 
The set of stochastic matrices is now replaced by the unitary group
$$\bU(n):=\{ A\in Mat_{n\times n}(\BC) \mid \ |Av|=|v|\ \forall \bra v\in \BC^n\}.
$$

The unitary group satisfies the third wish on the list:
   For all $A\in \bU(n)$, there exists a matrix $B\in \bU(n)$ satisfying
$B^2=A$.

Consider  wish 1:
  it is an open question! However   at least  our generalized
probabilistic computation includes our old probabilistic computation because
  $H$ and the    matrices \eqref{tofmat}, \eqref{negationmat}  are unitary. 

   \medskip

An indication that   generalized probability may be related
to quantum mechanics is that  the  interference patterns observed in the famous two slit 
experiments is manifested in  generalized probability:
one obtains a \lq\lq random bit\rq\rq\  
by applying $H$ to $\bra 0$: $H\bra 0=\frac 1{\sqrt 2}(\bra 0 + \bra 1)$.
However, if one applies a second quantum coin flip, one looses the randomness as
$H^2\bra 0=\bra 1$, which, as pointed out in \cite{MR3058839}, could be interpreted as a manifestation
of interference.

\subsection{Postulates of quantum mechanics and relevant linear algebra}

Here are the standard postulates of quantum mechanics and relevant definitions from linear algebra.

\noindent{\bf P1}. Associated to any isolated physical system is    a Hilbert space  $\cH$, called the {\it state space}.
The system is completely described at a given moment by a  unit vector $\bra\psi\in \cH$, called its {\it state vector}, which is well defined up to a phase $e^{i\th}$ with $\th\in \BR$. 
Alternatively one may work in projective space $\BP \cH$. 

\medskip

\noindent{\bf Explanations}. A  {\it Hilbert space} $\cH$ is a (complete) complex vector space endowed with a non-degenerate  Hermitian inner-product,
$h: \cH\times \cH\ra \BC$, where by definition $h$ is linear in the first factor
and conjugate linear in the second, $h(\bra v,\bra w)=\ol{h(\bra w,\bra v)}$ for all $v,w$, and $h(\bra v,\bra v)>0$ for
all $\bra v\neq 0$. 

\medskip

\noindent{\bf Remarks}:

When studying quantum mechanics in general, one needs to
allow infinite dimensional Hilbert spaces, but in the case of quantum computing, one restricts to finite
dimensional Hilbert spaces, usually $(\BC^2)^{\ot N}$. 

The Hermitian inner-product  $h$ allows an identification of $\cH$ with $\cH^*$ by $\bra w \mapsto \ket w:=h(\cdot, \bra w)$.
This identification will be used repeatedly.  
Write $h(\bra v,\bra w)=\braket wv$ and  $|v|=\sqrt{\braket v v}$ for the {\it length} of $\bra v$.

If $\cH= \BC^n$   with its standard basis, where 
$\bra v=(v_1\hd v_n) $,   the {\it standard Hermitian inner-product} on 
$\BC^n$ is $\braket w v=\sum_{j=1}^n \ol{w}_j v_j$. I will always assume $\BC^n$ is equipped with its
standard Hermitian inner-product. 

\medskip

\noindent{\bf P2}. The state of
an isolated  system   evolves with time according to the {\it Schr\"odinger equation}
$$
i\hbar \frac{d\bra\psi}{dt} = X\bra\psi
$$
where $\hbar$ is a constant ({\it Planck's constant}) and $X$ is a fixed {\it Hermitian operator},
called the {\it Hamiltonian} of the system.
Here, recall that  the {\it adjoint} of an operator $X\in \tend(\cH)$, is
 the operator $X^\dag\in \tend(\cH)$
such that $\braket{X^\dag v}{w}=\braket{v}{Xw}$ for all $v,w\in \cH$ and  $X$ is {\it Hermitian}
if $X=X^\dag$.
For a general Hilbert space, 
  the   Unitary group  is
$\bU(\cH):=\{ U\in \tend(\cH)\mid |Uv|=|v|\  \forall \bra v  \in \cH\}$.

How is   generalized probability related   to Schr\"odinger's equation?
Let  $U(t)\subset \bU(\cH)$ be a smooth curve  with $U(0)=\Id$.
Write $U'(0)=\frac d{dt}|_{t=0} U(t)$. 
Consider
\begin{align*}
0&= \frac d{dt}|_{t=0}\braket vw\\
&= \frac d{dt}|_{t=0}\braket  {U(t)v}{U(t)w}\\
&= \braket{U'(0)v} w + \braket v{U'(0)w}.
\end{align*}

Thus $iU'(0)$  is Hermitian. 
We are almost at Schr\"odinger's equation.  Let $\fu(\cH)=T_{\Id}\bU(\cH)$ denote the Lie
algebra of $\bU(\cH)$ so  $i \fu(\cH)$ is the space of Hermitian endomorphisms.
For $X\in \tend(\cH)$, write $X^k\in \tend(\cH)$ for $X\cdots X$ applied $k$ times.
Write $e^X:=\sum_{k=0}^\infty \frac 1{k!} X^k$.   If $X$ is Hermitian, then $e^{iX}\in \bU(\cH)$. 
Postulate 2 implies the system will evolve unitarily, by (assuming one starts at $t=0$),
$\bra{\psi_t}=U(t)\bra{\psi_0}
$, where 
$$
U(t)=e^{\frac{-itX }{\hbar}}.
$$

\subsection*{Measurements}
Our first two postulates  dealt  with isolated systems. In reality, no system is isolated and the whole
universe is modeled by one enormous Hilbert space. In practice, parts of the system  
are sufficiently isolated that they can be treated as isolated systems. However, they are occasionally
acted upon by the outside world, and one needs a way to describe this outside interference.
For our purposes, the isolated systems will be the Hilbert space attached to the input in a quantum algorithm
and the outside interference will be the measurement at the end.
That is, after a sequence of unitary operations one obtains a vector $\bra \psi=\sum z_j \bra j$ (here implicitly assuming the
Hilbert space is of countable dimension), and
as in generalized probability: 

{\bf P3}
The probability of obtaining outcome $j$ under a measurement equals $|z_j|^2$.

\medskip

In \S\ref{quantumreform}, motivated again by probability,   {\bf P1,P3} 
will be generalized to new postulates
that give rise to the same theory, but are more convenient to work with in information theory.   

\medskip

A typical situation in quantum mechanics and quantum computing is that there are
two or more isolated systems, say $\cH_A,\cH_B$ that are brought together (i.e., allowed to interact with each
other) to form a larger isolated system $\cH_{AB}$. The larger system is called the {\it composite system}.
In classical probability, the composite space is $\{ 0,1\}^{N_A}\times \{ 0,1\}^{N_B}$.  
In our generalized probability, the   composite space is $(\BC^2)^{\ot N_A}\ot (\BC^2)^{\ot N_B}=(\BC^2)^{\ot (N_A+N_B)}$:

{\bf P4}:  The state of a composite system $\cH_{AB}$ is the tensor product of the state
spaces of the component physical systems $\cH_A,\cH_B$: $\cH_{AB}=\cH_A\ot \cH_B$.

When dealing with composite systems, we will need to  allow partial measurements whose outcomes
are of the form $\bra I\ot \bra\phi$ with $\bra\phi$ arbitrary.

This tensor product structure gives rise to the notion of {\it entanglement}, which   
  accounts for phenomenon outside of our classical intuition, as discussed in the next section.

\begin{definition}
A state $\bra\psi\in \cH_1\otc \cH_n$ is called {\it separable} if it corresponds to a rank one tensor,
i.e., $\bra\psi=\bra {v_1}\otc \bra {v_n}$ with each $\bra{v_j}\in \cH_j$. Otherwise it
is {\it entangled}.
\end{definition}

\section{Entanglement phenomena}\label{EP}

\subsection{Super-dense coding}\footnote{Physicists use the word \lq\lq super\rq\rq\ in
the same way American teenagers use the word \lq\lq like\rq\rq .}
Physicists    describe their experiments in terms of two characters, Alice and Bob.
I generally  follow this convention.
Let $\cH=\BC^2\ot \BC^2=\cH_A\ot \cH_B$, and let $\epr= \frac{\bra{00}+\bra{11}}{\sqrt 2}$
(called  the {\it EPR state} in the physics literature after Einstein-Podolsky-Rosen).
Assume this state has been created,  both Alice and Bob are aware of it,  
   Alice is in possession of  the first qubit, and Bob the second. In particular Alice can act on the first
   qubit by unitary matrices and Bob can act on the second.
This all happens before the experiment begins.

Now say Alice wants to transmit a two  classical bit message to Bob, i.e., one of
the four states $\bra{00},\bra{01},\bra{10},\bra{11}$  by transmitting qubits.
We will see that she can do so transmitting just one qubit. 
If she manipulates her qubit by acting on the first $\BC^2$ by
a unitary transformation, $\epr$ will be manipulated.
She uses the following matrices depending on the message she wants to transmit:
$$
\begin{array}{|r|c|c|c|}
\hline
{\rm to \ transmit} & {\rm\ act\ by} & {\rm to \ obtain} \\
\hline
\bra{00}& \Id & \frac{\bra{00}+\bra{11}}{\sqrt 2}\\
 \bra{01}& \begin{pmatrix} 1&0\\ 0 & -1\end{pmatrix}=:\s_z & \frac{\bra{00}-\bra{11}}{\sqrt 2}\\ 
\bra{10} &   
\begin{pmatrix} 0&1\\ 1 & 0\end{pmatrix}=:\s_x & \frac{\bra{10}+\bra{01}}{\sqrt 2}\\
  \bra{11}&   
\begin{pmatrix}  0 & -1\\ 1 & 0\end{pmatrix}=:-i\s_y
& \frac{\bra{01}-\bra{10}}{\sqrt 2}
\end{array}
$$
where the names $\s_x,\s_y,\s_z$ are traditional in the physics literature (the {\it Pauli matrices}).
If Alice sends Bob her qubit, so he is now in possession of the modified
$\epr$ (although he does not see it), he can determine which of the four messages she sent him by measuring
$\epr$.
More precisely, first Bob acts on $\BC^2\ot \BC^2$ by a unitary transformation that takes the orthonormal
basis in the \lq\lq to obtain\rq\rq\ column to the standard orthonormal basis (this is a composition of the inverses
of two Hadamard matrices), to obtain a state
vector whose probability is concentrated at one of the four classical states, then he measures, and
obtains
 the correct classical state with probability one. 

In summary,  with preparation of an EPR state in advance, plus transmission of a single
qubit, one can transmit two classical bits of information.

\subsection{Quantum teleportation}\label{teleporationsect} A similar phenomenon is the misleadingly named 
{\it quantum teleportation},
where again Alice and Bob share half of an EPR state,
Alice is in possession of a qubit $\bra\psi=\a\bra 0 +\b\bra 1$, and wants to \lq\lq send\rq\rq\ $\bra\psi$ to Bob. However 
Alice is only allowed to transmit classical information to Bob (in particular, that information is transmitted
at a speed slower than the speed of light, which is why the use of the word
\lq\lq teleportation\rq\rq\ is misleading).  We will see that she can transmit a qubit to Bob by   transmitting
two classical bits.
Write the state of the system as
$$
\frac 1{\sqrt 2}\left[
\a\bra{0}\ot (\bra{00}+\bra{11})+\b\bra{1}\ot (\bra{00}+\bra{11})
\right]
$$
where Alice can operate on the first two qubits.
If Alice     acts on the first two qubits by   
$H\ot \s_x= \frac 1{\sqrt 2}\begin{pmatrix} 1&1\\ 1&-1\end{pmatrix}\ot \begin{pmatrix} 0&1\\ 1&0\end{pmatrix}$,
she obtains
$$
\frac 12\left[
\bra{00}\ot (\a\bra 0+\b\bra 1) +\bra{01}\ot(\a\bra 1+\b\bra 0)
+\bra{10}\ot (\a\bra 0-\b\bra 1)+ \bra{11}\ot (\a\bra 1-\b\bra 0)
\right].
$$
 
Notice that Bob's coefficient of Alice's $\bra{00}$ is the state $\bra\psi$ that is to be
transmitted.   Alice performs a measurement. 
If she has the good luck to obtain $\bra{00}$,
then she knows Bob has $\bra\psi$ and she can tell him classically that
he is in possession of $\bra\psi$. But say she obtains the state
$\bra{01}$: the situation is still good, she knows Bob is in possession
of a state such that,  if he acts on it with $\s_x=\begin{pmatrix}0&1\\ 1& 0\end{pmatrix}$,
he will obtain the state $\bra\psi$, so she just needs to tell him classically to apply $\s_x$.
Since they had communicated the algorithm in the past, all Alice really needs
to tell Bob in the first case is the classical message $00$ and in the second case the
message $01$. The   cases of  $10$ and $11$ are similar.

In summary, a shared EPR pair plus sending two classical bits of information
allows transmission of  one qubit. 
  
\subsection{Bell's game}
The 1935 Einstein-Podolsky-Rosen paper \cite{EPR} challenged quantum mechanics with the following thought experiment
that they believed implied instantaneous communication across distances, in violation of principles of relativity:
Alice and Bob prepare   $\epr=\eprx$, then travel far apart. Alice measures her bit. If she gets $0$,
then she can predict with certainty that Bob will get $0$ in his measurement, even if his measurement is taken
a second later and they are a light year apart.  

Ironically, this thought experiment has been made into an actual experiment.
One  modern interpretation (see, e.g., \cite{MR2500087}) is that there is no paradox because the system does not transmit
information faster than the speed of light, but rather they are acting on information that
has already been shared. What follows is a
version from 
\cite{PhysRevLett.23.880}, adapted from the presentation in \cite{MR2500087}.

Charlie chooses $x,y\in\{ 0,1\}$ at random and sends $x$ to Alice and $y$ to Bob.
Based on this information, Alice and Bob, without communicating with
each other,  get to choose bits $a,b$ and send them to
Charlie. The game is such that Alice and Bob play on a team. They win if $a\oplus b=x\ww y$, i.e., either $(x,y)\neq (1,1)$ and $a=b$
or $(x,y)=(1,1)$ and $a\neq b$.

\subsubsection{Classical version}

Note that if Alice and Bob both always choose $0$, they win with probability $\frac 34$.

\begin{theorem}\cite{Bell} Regardless of the strategy Alice and Bob use, they never 
win with probability greater than $\frac 34$.
\end{theorem}

\subsubsection{Quantum version} Although there is still no communication allowed between Alice and Bob, they will exploit
a pre-shared $\epr$ to gain an advantage over the classical case.
Alice and Bob prepare $\epr= \frac{\bra{00}+\bra{11}}{\sqrt 2}$ in advance, and  
Alice takes the first qubit and Bob the second.
When Alice gets $x$ from Charlie, if
$x=1$, she applies a rotation by $\frac{\pi}8$ to her qubit,   and if $x=0$ she does nothing.
When Bob gets $y$ from Charlie, 
 he applies a rotation by $-\frac{\pi}8$ to his qubit if
$y=1$ and if $y=0$  he does nothing. (The order these rotations are applied does not matter because
the corresponding operators on $(\BC^2)^{\ot 2}$ commute.)
Both of them measure their respective qubits and send the values obtained to Charlie.

\begin{theorem} With this strategy, Alice and Bob win with probability at least $\frac 45$.
\end{theorem}

The idea behind the strategy is that when $(x,y)\neq (1,1)$, the states of the two qubits
will have an angle at most  $\frac \pi 8$ between them, but when $(x,y)=(1,1)$, the angle
will be $\frac\pi 4$. 
 
\section{Quantum algorithms}\label{Qalgs}

\subsection{Grover's search algorithm}
The problem: given $F_n: \BF^n_2\ra \BF_2$, computable by a $poly(n)$-size classical circuit, 
find $a$ such that $F_n(a)=1$ if such $a$ exists.

Grover found a quantum circuit of size $poly(n)2^{\frac n2}$ that solves this problem with high probability.
Compare this with a brute force search, which requires a circuit of size $poly(n)2^n$. No classical or
probabilistic algorithm is known that does better than 
$poly(n)2^n$. Note that it also gives a size $poly(n)2^{\frac n2}$ probabilistic solution to the $\np$-complete
problem SAT (it is stronger, as it not only determines existence of a solution, but finds it).

I   present the algorithm for the following simplified version where one is promised
there exists exactly one solution. All essential ideas of the general case are here.

Problem: given $F_n: \BF^n_2\ra \BF_2$, computable by a $poly(n)$-size classical circuit, 
and the information that there is   exactly one vector $a$ with $F_n(a)=1$,  find $a$.

The idea of the algorithm is to start with a vector equidistant from all possible solutions,
and then to incrementally rotate it towards $a$.
What is strange for our classical intuition is that one is  able to rotate towards
the solution without knowing what it is, and similarly, we won't \lq\lq see\rq\rq\ the rotation
matrix either. 

Work in $(\BC^2)^{\ot n +s}$ where $s=s(n)$ is the size of the classical circuit needed to compute $F_n$.
I suppress reference to the $s$ \lq\lq workspace bits\rq\rq\ in what follows.
 
The following vector   is the average of all the classical (observable) states:
\be\label{avgclas}
\bra\av:=\frac 1{2^{\frac n2}}\sum_{I\in \{ 0,1\}^n } \bra I.
\ene

  To prepare $\bra{\av}$,    note that $H\bra 0=\frac 1{\sqrt 2}(\bra 0+\bra 1)$, so applying $H^{\ot n}$
to $\bra{0\cdots 0}$ transforms it to $\bra \av$. 
The cost of this is $n$ gates (matrices).

Since $\bra\av$  is   equidistant from all possible solution vectors,     $\braket\av a=\frac 1{2^{\frac n2}}$.
We want to rotate $\bra\av$ towards the unknown $a$. 
Recall that $\tcos(\angle (\bra v,\bra w))=\frac{\braket vw}{|v||w|}$. 
Write the angle between $\av$ and $a$  as $\frac \pi 2 - \th$, so $\tsin(\th)=\frac 1{2^{\frac n2}}$.

A rotation is a product of two reflections.
  In order to perform the rotation $R$ that moves $\bra\av$ towards $\bra a$,
 first reflect in the hyperplane orthogonal to $\bra a$, and then in the hyperplane orthogonal to $\bra\av$.

Consider the map    
\be\label{clasquan}\bra {xy}\mapsto \bra{x(y\op F(x))}
\ene  
defined on basis vectors and extended linearly. 
 To execute this,   use the $s$ workspace bits
that are suppressed from the notation, to effect $s$ reversible classical gates.   Initially
set $y=0$  so that  the image is $\bra{x0}$ for $x\neq a$,
and $\bra{x1}$ when $x=a$. 
Next apply the  quantum gate $\Id\ot \begin{pmatrix} 1&0\\ 0&-1\end{pmatrix}$
  which sends $\bra{x0}\mapsto \bra{x0}$,
and $\bra{x1}\mapsto -\bra{x1}$. 
Finally apply the map $\bra {xy}\mapsto \bra{x(y\op F(x))}$ again.

Thus $\bra{a0}\mapsto -\bra{a0}$ and all other vectors $\bra{b0}$  are mapped to themselves, as desired.

Next we need to  reflect around $\bra\av$.
It is easy to reflect around a classical state, so first
  perform the map
$H^{-1\ot n}$   
that sends $\bra \av$  to $\bra{0\cdots 0}$, then reflect 
in the hyperplane perpendicular to $\bra{0\cdots 0}$ using the 
Boolean function $g: \BF_2^n\ra \BF_2$ that
outputs $1$ if and only if its input is $(0\hd 0)$, in the role of $F$ for our previous reflection, then apply Hadamard again
so the resulting reflection is about $\bra\av$.

The composition of these two reflections is the desired rotation $R$. The vector $R\bra \av$ is not useful as measuring it
only slightly increases the probability of obtaining $\bra a$, but if one
composes  $R$ with itself $O(\frac 1{\th})$ times,
  one obtains a vector much closer to $\bra a$. (Note that $\th\sim \tsin(\th)$ so $\frac 1\th\sim \sqrt{N}$.)

\subsection{The quantum discrete Fourier transform}
Underlying the famous quantum algorithm of Shor for factoring integers 
and Simon's algorithm that led up to it, are \lq\lq quantum\rq\rq\ versions of the discrete Fourier transform on finite
abelian groups.  

The DFT for $ \BZ/M\BZ$,   in vector notation,
for $j\in \BZ/M\BZ$, is
$$
\bra j\mapsto  \frac 1{\sqrt{M}} \sum_{k=0}^{M-1}
\o^{jk}\bra k
$$
where $\o=e^{\frac{2\pi i}M}$.
It is a unitary change of basis such that in the new basis,
multiplication in $\BZ/M\BZ$ is given by a diagonal matrix, and the classical
FFT writes the DFT as a product of $O(\tlog(M))$ sparse matrices (each with $M<< M^2$ nonzero entries),
for a total cost of $O(\tlog(M)M)<O(M^2)$ arithmetic operations.

Write $M=2^m$.
The DFT  can be written as a product of $O(m^3)=O(\tlog(M)^3)$   controlled local  unitary
operators.
Hence one can approximately obtain the output vector   by a sequence of  $poly(m)$ 
unitary operators from our gate set with the caveat that we won't be able to
\lq\lq see\rq\rq\ it. 

Here is the quantum DFT:
It will be convenient
to express $j$ in binary and view $\BC^M=(\BC^2)^{\ot m}$, i.e.,
write
$$
\bra j=\bra{j_1}\otc \bra{j_m}
$$
where $j=j_12^{m-1}+j_2 2^{m-2}+\cdots + j_m 2^0$ and $j_i\in \{0,1\}$.
Write the DFT as
\begin{align}
\nonumber &\bra{j_1}\otc \bra{j_m} \\
\nonumber & \mapsto \frac 1{\sqrt{M}} \sum_{k=0}^{M-1}\o^{jk}\bra k\\
\nonumber &=\frac 1{\sqrt{M}} \sum_{k_i\in \{0,1\} }\o^{j(\sum_{l=1}^m k_l2^{m-l})}\bra {k_1}\otc \bra{k_m}\\
\nonumber &=\frac 1{\sqrt{M}} \sum_{k_i\in \{0,1\} }
\bigotimes_{l=1}^m \left[ \o^{j  k_l2^{m-l}}\bra {k_l}\right] \\
\nonumber &=\frac 1{\sqrt{M}} \sum_{k_i\in \{0,1\} }
\bigotimes_{l=1}^m \left[ \o^{(j_12^{2m-1-l}+\cdots + j_m2^{m-l})  k_l }\bra {k_l}\right] \\
\label{quandft}&=\frac 1{2^{\frac m2}}
(\bra 0 + \o^{  {j_m}2^{-1}}\bra 1)
\ot
(\bra 0 + \o^{  {j_{m-1}}2^{-1}+  {j_m}2^{-2}}\bra 1)
\ot 
(\bra 0 + \o^{  {j_{m-2}}2^{-1} +  {j_{m-1}}2^{-2}   + {j_{m }}2^{-3}}\bra 1)\\
&\nonumber\ \ \otc 
(\bra 0 + \o^{\sum_{s=0}^{m-1}  {j_{m-s}} 2^{m-(s+1)}   } \bra 1)
\end{align}
where for the last line if $2m-s-l>m$, i.e., $s+l<m$, there is no contribution  with $j_s$ because
  $\o^{2^m}=1$, and I multiplied all terms by $1=\o^{2^{-m}}$ to have negative exponents.
  
It will be notationally more convenient to write the quantum circuit for this vector with the order of factors reversed, so  I describe a quantum circuit that produces
\begin{align}
\label{twisteddft}
&\frac 1{\sqrt{2}}(\bra 0 + \o^{\sum_{s=0}^{m-1}  {j_{m-s}} 2^{m-(s+1)}   } \bra 1)  
  \otc\frac 1{\sqrt{2}}(\bra 0 + \o^{  {j_{m-2}}2^{-1} +  {j_{m-1}}2^{-2}   + {j_{m }}2^{-3}}\bra 1)\\
&\nonumber \ot
\frac 1{\sqrt{2}}(\bra 0 + \o^{  {j_{m-1}}2^{-1}+  {j_m}2^{-2}}\bra 1)
\ot 
\frac 1{\sqrt{2}} (\bra 0 + \o^{  {j_m}2^{-1}}\bra 1).
\end{align}
  
Set 
\be\label{rk}
R_k=\begin{pmatrix} 1&0\\ 0 & \o^{2^k}\end{pmatrix},
\ene
then \eqref{twisteddft}  is obtained  as follows:
first apply $H$ to $(\BC^2)_1$ then  a linear map
  $\La 1 R_j$,  defined by $\bra x\ot \bra y\mapsto \bra x\ot R_j\bra y$ if $\bra x\neq \bra 0$ and
 to $\bra x\ot \bra y$ if $\bra x=\bra 0$, to $(\BC^2)_j\ot (\BC^2)_1$ for $j=2\hd m$.
 Note that at this point only the $(\BC^2)_1$-term has been altered.
 From now on  leave the $(\BC^2)_1$-slot alone.
Next apply $H$ to $(\BC^2)_{2}$ then   $\La 1 R_{j-1}$ to $(\BC^2)_{j}\ot (\BC^2)_{2}$ for $j=3\hd m$.
Then 
apply $H$ to $(\BC^2)_{3}$ then   $\La 1 R_{j-2}$ to $(\BC^2)_j\ot (\BC^2)_3$ for $j=4\hd m$.
Continue, until finally one just applies $H$ to $(\BC^2)_m$.
Finally to obtain the DFT, reverse the orders of the factors (a classical operation).

In practice, one has to fix a quantum gate set, i.e., a finite set of unitary
operators that will be allowed in algorithms,  in advance, so in general it will
be necessary  to approximate
the transformations $R_k$ from elements of our gate set, so one only obtains
an approximation of  the DFT. 

\subsection{The hidden subgroup problem} 
 Given  a discrete group $G$ with
a specific representation of its elements in binary,
a  function $f: G\ra \BF_2^n$, and a device that computes $f$ (for unit cost),  and the knowledge that there exists
a subgroup $G'\subset G$ such that $f(x)=f(y)$ if and only if
$xy\inv \in G'$, find $G'$.

For finitely generated abelian groups, it is sufficient to solve the problem for $G=\BZ^{\op k}$ as all 
finitely generated abelian groups are quotients of some $\BZ^{\op  k}$.

Simons algorithm is for the hidden subgroup problem with   $G=\BZ_2^{\op m}$.
The $DFT_2$ matrix is just
$$
H=\frac 1{\sqrt 2}
\begin{pmatrix}
1& -1\\ -1 & 1\end{pmatrix}
$$
and $G'$ is the subgroup generated by $a\in \BZ_2^{\op m}$.

Shor's algorithm for factoring (after classical preparation) amounts to  the case $G=\BZ$ and $F$ is the function $x\mapsto a^x\tmod N$.

\section{Classical information theory}\label{cinfo}
Quantum information theory is based on classical information theory, so I review the
classical theory. 
  The discovery/invention of the bit by Tukey and its
development by Shannon \cite{MR0026286}
was one of the great scientific
achievements of the twentieth century, as it changed the way one views information,
giving it an abstract formalism that is discussed in this section.
The link to quantum information is explained in \S\ref{quantumchan}.

\subsection{Data compression: noiseless channels}  
  (Following \cite{2016arXiv160401790B})
A source emits symbols $x$ from an alphabet $\cX$ that
we want to store efficiently
so we try to encode $x$ in a small number of bits, to say $y\in \cY$ in a way that
one can decode it later to recover $x$.

\begin{figure}[!htb]\begin{center}
\includegraphics[scale=.3]{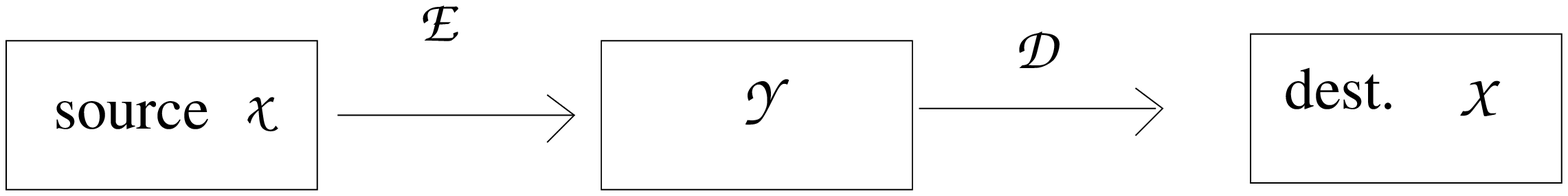}
\caption{\small{Message from source encoded into bits then decoded}}  \label{fig:encodedecode}
\end{center}
\end{figure}

The symbols from $\cX$ do not necessarily occur with the same frequency.
Let $p=P_{\cX}$ denote the associated probability distribution. What is
  the minimum possible size for $\cY$? Since we are dealing in bits, it will
be convenient to use the   logarithms of cardinalities, so define the {\it capacity} as 
$\capp(P_{\cX}):=\tmin\tlog|\cY|$.

Consider the case $\cX=\{a,b,c,d\}$ where $p(a)=0.1$, $p(b)=0$,
$p(c)=0.4$ and $p(d)=0.5$. One  can clearly get away with $|\cY|=3$,
  e.g., for the encoder, send $a,b$ to $1$,
$c$ to $2$ and $d$ to $3$, then for the decoder, send
$1$ to $a$, $2$ to $c$ and $3$ to $d$. 
In general, one can always throw away symbols with probability zero.
On the other hand, one cannot map two distinct symbols that do occur to
the same symbol, as there would be no way to distinguish them when decoding.
Thus   $\capp(p)= \tlog\tsupp(p)$,
where $\tsupp(p)=\# \{ x\in \cX \mid p(x)>0\}$.

Now say we are willing to tolerate a small error. First 
rephrase what we did probabilistically:
Let $p^{enc}(y|x)$ denote the conditional probability distribution of
the encoder $\cE$ and $p^{dec}(x|y)$ that of the decoder $\cD$.
Our requirement was for all $x$,
$$
p[x=\cD\circ \cE(x)]=\sum_{ y,x'} 
 p^{enc}(y|x)p^{dec}(x'|y)\d_{x,x'}=1.
$$
Now relax it to 
$$
 \sum_{x,y,x'} 
p(x)p^{enc}(y|x)p^{dec}(x'|y)\d_{x,x'}\geq 1-\ep.
$$
for some error $\ep$ that we are willing to tolerate.
In addition to throwing out the symbols that do not appear, we may also
discard the largest set of symbols whose total probability is smaller than $\ep$.
Call the corresponding quantity $\capp^{\ep}(p)$.

In the example above, if one takes $\ep> 0.1$, one can lower storage cost,
taking $|\cY|=2$.

Recall that a probability distribution $p$ on $\cX$ must satisfy $\sum_{x\in \cX}p(x)=1$.
Relax this to {\it non-normalized} probability distributions, $q$, where
$q(x)\geq 0$ for all $x\in \cX$ and $\sum_{x\in \cX}q(x)\leq 1$.
We obtain: 
$\capp^{\ep}(p)= \tmin \tlog\tsupp(q)$, where the min is taken over all
non-normalized probability distributions $q$ satisfying $q(x)\leq p(x)$ and  $\sum_{x\in \cX} q(x)\geq 1-\ep$.

Now say we get not a single symbol, but a string of $n$ symbols, so we seek
an encoder $\cE: \cX^n\ra \cY(n)$, where $\cY(n)$ is a set that
varies with $n$,   and decoder $\cD: \cY(n)\ra \cX^n$, and we want to 
minimize $|\cY(n)|$, with a tolerance of error that goes to zero as $n$ goes  
to infinity.   In practice one wants to send
information 
  through a communication channel (e.g. telegraph wire). 
The channel can only send a limited number of bits per second, and we
want to maximize the amount of information we can send per second.
Define  $\trate(p):=\tlim_{\ep\ra 0}\tlim_{n\ra \infty} \frac 1n \capp^{\ep}(p^n)$.  

The string $x_1\cdots x_n=:\ol x^n$ is identically and independently distributed (i.i.d),
that is each $x_j$ is drawn from the same probability distribution and the draw of
$x_j$ is independent of the draws of the other $x_i$.
Say $\cX=\{1\hd d\}$ with $p(j)=p_j$. The probability of any given
string occurring depends only on the number of $1$'s $2$'s etc..  in the
string and not on
their order. A string with $c_j$ $j$'s occurs with probability
$p_1^{c_1}\cdots p_d^{c_d}$. (Note that $c_1+\cdots + c_d=n$.)  The number of strings with this probability is
$$
\binom n {c_1\hd c_d}:=\frac{n!}{c_1!\cdots c_d!}
$$
and we need to estimate this quantity.

Stirling's formula implies
$
\tlog(n!)=n\tlog(n)-n+ O(\tlog(n))$. 
In particular, for $0<\b<1$ such that $\b n\in \BZ$,
$$
\tlog\binom n{\b n}
 =
n[-\b \tlog(\b) -(1-\b)\tlog(1-\b)]+O(\tlog(n)).
$$
Let $H(\b)=-\b \tlog(\b) -(1-\b)\tlog(1-\b)$ and more
generally, for $\ol{p}=(p_1\hd p_d)$, let
$$H(\ol{p})=-\sum_{i=1}^d p_i\tlog(p_i),
$$
the {\it Shannon entropy} of $\ol{p}$.
It   plays a central role in information theory.

Define a map $wt: \cX^n\ra \BR^d$ by $\ol x^n\mapsto (c_1\hd c_d)$, where $c_j$ is the number of 
$j$'s appearing in $\ol x^n$. 
Then  the expectation is $E[wt(\ol x^n)]=(np_1\hd np_d)$. 
The weak law of large numbers   states that for any $\ep>0$,  
$$
\tlim_{n\ra \infty} p[|| \frac 1n  wt(\ol{x}^n)-E[wt(\ol x^n)) ]||_1>\ep]=0
$$
where for  $f: \cZ\ra   \BR^d$, define  $||f||_1=\sum_{z\in \cZ} |f(z)|$. In our case, $\cZ=\cX^n$.

Now simply throw out all strings $\ol{x}^n$ with
$||\frac{1}n (wt(  \ol x^n)-E[wt(\ol   x^n)) ]||_1>\ep$, and  take $\cY(n)$ of size
\begin{align*}
|\cY(n)|&=\#\{ \ol x^n\mid ||\frac{1}n (wt(\ol x^n)-E[wt(\ol x^n))]||_1<\ep \}\\
&=
\sum_{\stack{\ol x^n\mid}{ ||\frac{1}n( wt(\ol x^n)-E[wt(\ol x^n))]||_1<\ep }} \binom{n}{wt(\ol{x}^n) }.
\end{align*}

If $\ep$ is small, the multinomial coefficients appearing will all be very close
to
$$
\binom n{np_1\hd np_d}
$$
and for what follows, one can take the crude approximation \be\label{ysize}
|\cY(n)|\leq  poly(n) \binom n{np_1\hd np_d}
\ene
(recall that  $d$ is fixed).

Taking 
logarithms, the right hand side of \eqref{ysize} becomes
$nH(\ol{p})+O(\tlog(n))
$.
Thus
$$
\frac 1 n \tlog |\cY(n)|\leq H( \ol{p})+o(1)
$$
and $\trate(\ol{p})\leq H(\ol{p})$.

\begin{theorem} \cite{MR0026286} $\trate(\ol{p})= H(\ol{p})$.
\end{theorem}

The full proof uses the strong law of large numbers.

\subsection{Transmission over noisy channels}

Say   symbols $x$ are transmitted over a channel subject to
noise, and   symbols $y$ are received  so one may or may not have $y=x$. Intuitively, if the noise is small, with some redundancy
it should be possible  to communicate accurate messages most of the time.
Let $\trate$ denote the maximal possible rate of transmission. In a noiseless channel
this is just $H(p_{\cX})$, but now we must subtract off something to account for the uncertainty
that, upon receiving  $y$, that it was the signal   sent. This something
will be the {\it conditional entropy}:
  Recall the conditional probability of $i$ occurring given
knowledge that $j$ occurs (assuming $p(j)>0$):  $p_{\cX | \cY}(i|j)=\frac{p_{\cX, \cY}(i,j)}{p_{  \cY}(j)}$ 
(also recall $p_{  \cY}(j)=\sum_ip_{\cX, \cY}(i,j)$).
Define the conditional entropy
$$H(\ol p_\cY|\ol p_\cX ):=-\sum_{i,j}p_{\cX, \cY}(i,j)\tlog p_{\cY | \cX}(j|i).
$$
Note that
\be\label{conduncert}
H(\ol p_\cY|\ol p_\cX )=H(\ol p_{\cX,\cY} )-H(\ol p_{\cX } )
\ene
or equivalently  $H(\ol p_{\cX,\cY} )=H(\ol p_{\cX } )+H(\ol p_\cY|\ol p_\cX )$, the uncertainty of   $p_{\cX,\cY}$ is the uncertainty
of $p_{\cX }$ plus the uncertainty of $p_{\cY }$ given $p_{\cX }$.

In particular  $
H(\ol p_{\cY } )\geq H(\ol p_\cY|\ol p_\cX ) 
$, 
i.e., with extra knowledge, our uncertainty about $p_{\cY } $ cannot increase, and decreases unless
$p_{\cX } $ and $p_{\cY } $ are independent).

 \subsubsection{Capacity of a noisy channel}
 
 Define the {\it capacity} of a noisy channel to be the maximum rate over all possible probability distributions
 on the source:
 $$
 \capp:=\tmax_{q_{\cX}}\left(H(q_{\cX})-H(q_{\cX}|p_{\cY})\right).
 $$
  Shannon  \cite{MR0026286} proves that $\capp$ lives up to its name: if the entropy 
  of a discrete channel is below $\capp$
  then there exists an encoding $\ol p$ of the source such that information
 can be transmitted over the channel with an arbitrarily small frequency  of errors.  
  The basic idea is the same as the noiseless case, however there is a novel feature that
  now occurs frequently in complexity theory arguments - that instead of producing an
  algorithm to find the efficient encoding, Shannon showed that a {\it random} choice
  of encoding will work.   
  
 After presenting the proof, Shannon remarks:
\lq\lq An attempt to obtain a good approximation to ideal coding by following the method of the proof is
generally impractical. ... Probably this is no accident but is
related to the difficulty of giving an explicit construction for a good approximation to a random sequence\rq\rq .
To my knowledge, this is the first time   that the difficulty of \lq\lq finding hay in a haystack\rq\rq\
(phrase due to Howard Karloff) is mentioned in print. 
This problem is central to complexity: for example, Valiant's algebraic version of
$\p\neq\np$ can be phrased as the problem of finding a    sequence of
explicit polynomials that are difficult to compute, while it is known that
a random sequence is indeed difficult to compute.
According to A. Wigderson, the difficulty of writing down random objects
problem was also explicitly discussed by Erd\"os, in the context of random graphs, at least as early as 1947, in relation to his seminar paper
\cite{MR0019911}. This paper, along with \cite{MR0026286} gave rise to the
 now ubiquitous probabilistic method 
 in complexity theory.

\section{Reformulation of quantum mechanics}\label{quantumreform}
I discuss two inconveniences about our formulation of the postulates of quantum mechanics, leading
to a   formulation of the postulates in terms of density operators.

\subsection{Partial measurements} 
A measurement of a state
$\bra\psi =\sum z_I\bra I$ was defined  a procedure that gives us $I=(i_1\hd i_n)\in \{0,1\}^n$
with probability $|z_I|^2$. 
But in our algorithms, this is not   what happened: we  were working
not in $(\BC^2)^{\ot n}$, but $(\BC^2)^{\ot n+m}$ where there were $m$
\lq\lq workspace\rq\rq\ qubits   we were not interested in measuring. 
So our measurement was more like the projections onto the {\it spaces}
$\bra I\ot (\BC^2)^{\ot m}$.  I now define this generalized notion of
measurement.

To make the transition, first observe that 
$
  |z_I|^2 =\ket \psi \Pi_I \bra \psi$,
where $\Pi_I : (\BC^2)^{\ot n}\ra \BC\bra I$ is the orthogonal projection onto the line
spanned by $\bra I$.

Now say we are only interested in the first $n$ bits of a system of $n+m$ bits, and   want 
to know the probability a measurement gives rise to some  $I$ represented by a vector $\bra I\in (\BC^2)^{\ot n}$,
but we have $\bra \psi\in (\BC^2)^{\ot n+m}$. Adopt the notation $\ketbra \phi \psi:=\bra \phi\ot \ket\psi$.
Then the probability of obtaining $\bra I$ given $\bra \psi$ is
\begin{align}
\nonumber p(  \bra I \mid \bra \psi)&=\sum_{J\in\{0,1\}^m} p(\bra \psi, \bra{IJ})\\
\nonumber &= \sum_J \braket \psi{IJ}\braket{IJ}\psi\\
\nonumber &=\ket \psi (\bra I\ket I \ot \Id_{(\BC^2)^{\ot m}})\bra\psi\\
 \nonumber  &=\ket \psi \Pi_{\cM}\bra\psi
\end{align}
where $\Pi_{\cM}: (\BC^2)^{\ot n+m}\ra \bra I\ot (\BC^2)^{\ot m}=:\cM$ is
the orthogonal projection operator.
Then  $p(  \bra I \mid \bra \psi)=\ket \psi \Pi_{\cM}\bra\psi$.
With this definition, one can allow $\cM\subset \cH$ to be {\it any} linear
subspace, which will simplify our measurements. (Earlier,    if we wanted to measure
the probability of a non-basis state, we had to change bases before measuring.)

One may think of projection operators as representing outside interference of a quantum
system, like adding a filter to beams being sent that destroy states not in $\cM$.
 Recall that in classical probability, one has the identity:
\be
\label{totalprobcap}p(M_1\cup M_2)=
  p(M_1)+p(M_2) -
p(M_1\cap M_2)
.
\ene

The quantum analog is {\it false} in general:
Let $\cH=\BC^2$,
$\cM_1=\BC\bra 0$ and $\cM_2=\BC(\bra 0+\bra 1)$
Let $\bra\psi=\a\bra 0 +\b \bra 1$ with $|\a|^2+|\b|^2=1$.
Then (and in general) 
$p(\tspan\{ \cM_1,\cM_2 \})\neq 
 p( \cM_1)+
p( \cM_2 )  -p(  \cM_1\cap\cM_2 )$.

However, one  can recover \eqref{totalprobcap} if the projection operators commute:
\begin{proposition}\label{totalrecov}
 If $\Pi_{\cM_1}\Pi_{\cM_2}=\Pi_{\cM_2}\Pi_{\cM_1}$  then  
$p(\tspan\{ \cM_1,\cM_2 \})= 
 p( \cM_1)+
p( \cM_2 )  -p(  \cM_1\cap\cM_2 )$.
\end{proposition}

\subsection{Mixing classical and quantum probability}
A typical situation in probability is as follows: you want a cookie, but can't make up your
mind which kind, so you decide to take one at random from the cookie jar to eat. However when
 you open the cupboard, you find there are two different cookie jars $H$ and $T$, each with a different 
 distribution of cookies, say $P_H$ and $P_T$. You decide to flip a coin to decide which jar and
say your coin is biased with probability $p$ for heads (choice $H$).
The resulting probability distribution is
$$
pP_H +(1-p)P_T.
$$
Let's   encode this scenario with vectors.
Classically, if vectors corresponding to $P_H,P_T$ are respectively
$v_H,v_T$, the new vector is $pv_H+(1-p)v_T$.
The probability of drawing a chocolate chip (CC) cookie is
$pP_H(CC) +(1-p)P_T(CC)= pv_{H,CC}+(1-p)v_{T,CC}$.

But what should one take in generalized probability (where one uses  the $\ell_2$ norm
instead of the $\ell_1$ norm)? Given $\bra{\psi_A}=\sum z_I\bra I, \bra{\psi_B}=\sum w_J\bra J$,
we want to make a measurement that gives us $p |z_{CC}|^2+(1-p)|w_{CC}|^2$.
Unfortunately $|p z_{CC}+(1-p)w_{CC}|^2\neq  p |z_{CC}|^2+(1-p)|w_{CC}|^2$ in general.
To fix this problem I will enlarge the notion of state and further modify our notion of measurement.

Our problem comes from having a mixture of $\ell_1$ and $\ell_2$ norms. The  fix will be
to rewrite $\bra \psi$ in a way that the $\ell_2$ norm becomes an $\ell_1$ norm. 
That is,  I construct an object that naturally contains the squares of the norms of the coefficients
of $\bra {\psi_A}$. Consider the endomorphism $\ketbra{\psi_A}{\psi_A}=\sum_{I,J} z_I\ol{z}_J \ketbra IJ$.
It is rank one, and in our standard basis  its diagonal entries are the quantities we want.

To   measure  them, let $\Pi_J$ denote the projection onto the $J$-th coordinate.
Then
$$
\ttrace(\Pi_J \ketbra{\psi_A}{\psi_A})=|z_{A,J}|^2
$$
is the desired quantity.

Now back to our cookie jars, set
$$
\rho=p \ketbra{\psi_A}{\psi_A}+(1-p) \ketbra{\psi_B}{\psi_B}
$$
and observe that
$$
\ttrace(\Pi_J\rho)=  p |z_{A,J}|^2+(1-p)|z_{B,J}|^2 
$$
as desired.

Given
a finite set of states
$\{\bra{\psi_1}\hd \bra {\psi_s}\}$, with $p(\bra{\psi_i})=p_i$, and 
$\sum_ip_i=1$,
set $\rho=\sum_k p_k\ketbra{\psi_k} {\psi_k}\in \tend(\cH)$. Note that $\rho$ has the properties
\begin{enumerate}
\item $\rho=\rho^\dag$, i.e., $\rho$ is Hermitian,
\item $\forall \bra\eta$, $\braket{\eta|\rho} \eta\geq 0$, i.e., $\rho$ is {\it positive},
\item $\ttrace(\rho)=1$.
\end{enumerate}

 This motivates the
following definition:

\begin{definition}
An operator  $\rho\in \tend(\cH)$ satisfying 1,2,3 above is called a {\it density operator}.
\end{definition}

Note that a density operator that is diagonal in the standard basis of $\BC^d$
corresponds to a probability distribution on $\{1\hd d\}$, so the definition 
includes classical probability as well as our old notion of state (which are the rank one density operators).
The
set of density operators is invariant under   
the induced action of $\bU(\cH)$ on $\tend(\cH)$.

 Different scenarios can lead to the same density operator. However, 
two states with the same density operator are physically indistinguishable.

\subsection{Reformulation of the postulates of quantum mechanics} \ 

{\bf Postulate 1}.  Associated to any isolated physical system is a Hilbert space $\cH$,
call the {\it state space}. The system is described by its density operator
$\rho\in \tend(\cH)$. 

{\bf Postulate 2}.  The evolution of an isolated system is described by 
the action of unitary operators on $\rho$.

{\bf Postulate 3}. (POVM) Measurements correspond to a collection
of   projection operators $\Pi_{\cM_j}$ such that $\sum_k \Pi_{\cM_k}=\Id_{\cH}$.
The probability that $\rho$ is in measured in state $\cM_j$ is $\ttrace(\Pi_{\cM_j}\rho)$. 

Sometimes it is convenient to allow more general measurements:

{\bf Postulate 3'}. Projective measurements correspond to a collection
of  Hermitian   operators $X_j\in \tend\cH$ such that $\sum_k   X_k=\Id_{\cH}$.
The probability that $\rho$ is in measured in state $X_j$ is $\ttrace(    X_j\rho)$. 

Postulate 4 regarding composite systems is unchanged.

\begin{remark}
Note that for $A\in \tend\cH=\cH^*\ot \cH$, $\ttrace(A)$ is
the image of $A$ under the contraction map $\cH^*\ot \cH\ra \BC$,
$\ket v\ot \bra w\mapsto \braket v w$.
For $A\in \tend(\cH_1\ot \cH_2)=(\cH_1^*\ot \cH_2^*)\ot (\cH_1\ot \cH_2)$, define the partial trace
$\ttrace_{\cH_1}(A)$ to be the image of $A$ under the contraction
$\cH_1^*\ot \cH_2^*\ot \cH_1\ot \cH_2\ra \cH_2^*\ot \cH_2$ given by
$\ket \phi \ot \ket \psi \ot \bra v\ot \bra w\mapsto \braket \phi v\ket \psi \ot \bra w=\braket \phi v\ketbra w\psi $.
\end{remark}

\subsection{Expectation and the uncertainty principle}

Let $A\in \tend(H)$ be a Hermitian operator with eigenvalues $\l_1\hd \l_k$ and eigenspaces $\cM_j$. 
If our system is in state $\rho$, one  can consider $A$ as a random variable
that takes the value $\l_j$ with probability $\ttrace(\Pi_{\cM_j} \rho)$.

The expectation of a random variable $X: \cX\ra \BR$ is
$E[X]:=\sum_{j\in \cX} X(j)p(j)$.

If a system is in state $\rho$, the expectation   of   a Hermitian operator $A\in \tend(H)$ 
is $\ttrace(A\rho)$ because
$
 E[A]
 =\sum_{\l_j} \l_j \ttrace(\Pi_{\cM_j} \rho) 
=   \ttrace((\sum_{\l_j}\l_j\Pi_{\cM_j}) \rho) 
 =\ttrace(A\rho)$.

One way mathematicians describe the famous Heisenberg uncertainty principle is
that it is impossible to localize both a function and its Fourier transform.
Another interpretation comes from probability:

First note that given a random variable, or Hermitian operator $X$, one can replace it with
an operator of mean zero $\hat X:= X-E(X \rho)\Id$. For notational convenience, I state
the uncertainty principle for such shifted operators. 

 The variance $var(X)$ of a random variable is $var(X)=E[X-E(X)]^2$. 
The   standard deviation $\s(X)=\sqrt{var(X)}$ of   $X$ is
a measure of the failure of the corresponding probability distribution to be concentrated
at a point, i.e., failure of the induced probability distribution to have a certain outcome.

\begin{proposition}\label{heisprop}
Let $X,Y$ be Hermitian operators of mean zero, corresponding to observables on a system in state $\rho$, let
Then 
$$
\s(  X )\s(  Y )\geq \frac{|\ttrace   ([  X,  Y]\rho) |}2.
$$
\end{proposition}

The uncertainty principle says that the failure of two Hermitian operators
to commute lower bounds the product of their uncertainties. In particular,
if they do not commute, neither can give rise to a classical (certain) measurement.
It  is a consequence of the   Cauchy-Schwarz inequality.

\subsection{Pure and mixed states}
\begin{definition} Let $\rho\in \tend(\cH)$ be a density operator. If $\trank(\rho)=1$, i.e.
$\rho= \ketbra \xi\xi$, $\rho$ is called a {\it pure state}, and otherwise it is
called a {\it mixed state}.
\end{definition}

The partial trace of a pure state can be a mixed state. For
example, if $\rho=\ketbra \psi\psi$ with $\psi=\frac 1{\sqrt 2} (\bra {00} +\bra{11})\in \cH_1\ot \cH_2$,
then $\ttrace_{\cH_2}(\rho)=\frac 12(\ketbra 0  0 + \ketbra 11)$.

The following proposition shows that one could   avoid  density operators altogether by
working on a larger space:

\begin{proposition} \label{largespaceden} An arbitrary mixed state $\rho\in \tend(\cH)$ can be represented
as the partial trace $\ttrace_{\cH'}\ketbra \psi \psi $ of a pure state  in $\tend(\cH\ot \cH')$
for some Hilbert space $\cH'$.
In fact, one can always take $\cH'=\cH^*$.
\end{proposition}

Given a density operator $\rho\in \tend(\cH)$, there is
a well defined density operator $\sqrt{\rho}\in \tend(\cH)$ whose
eigenvectors are the same as for $\rho$, and whose eigenvalues
are the positive square roots of the eigenvalues of $\rho$.
To prove the proposition,    given $\rho\in \cH\ot \cH^*$,
consider $\ketbra {\sqrt \rho}  {\sqrt\rho}\in \tend(\cH\ot \cH^*)$.
Then $\rho= \ttrace_{\cH^*}(\ketbra {\sqrt \rho}  {\sqrt\rho})$.
 A pure state whose partial trace is $\rho$ is called a {\it purification} of $\rho$.

 \section{Communication across a  quantum channel}\label{quantumchan}
 Now instead of having a source $\cX^{\times n}$ our \lq\lq source\rq\rq\ is $\cH^{\ot n}$,
 where one can  think of $\cH^{\ot n}=\cH_A^{\ot n}$, and Alice will \lq\lq transmit\rq\rq\ a 
 state to Bob, 
 and instead of  a probability distribution $p$ one has a density operator $\rho$.

What is a quantum channel?
It should be a linear map sending
$\rho\in \tend(\cH_A)$ to some $\Phi(\rho)\in \tend(\cH_B)$.

First consider the special case $\cH_A=\cH_B$.
One should allow coupling with an auxiliary system, i.e.,
\be\label{cptp}
\rho\mapsto \rho\ot \s \in \tend(\cH_A\ot \cH_C).
\ene
One  should also allow the state $\rho\ot \s$ to evolve in $\tend(\cH_A\ot \cH_C)$, i.e., be acted upon by an
arbitrary $U\in \bU(\cH_A\ot \cH_C)$.
Finally one should allow measurements, i.e., tracing out the $\cH_C$ part.
In summary, a quantum channel $\cH_A\ra \cH_A$ is a map of the form $\rho\mapsto \ttrace_{\cH_C}(U(\rho\ot \s)U\inv)$. More generally to go from
$\cH_A$ to $\cH_B$, one needs to allow isometries as well.
Such maps are the  {\it completely positive trace preserving maps} (CPTP),
where a map $\Lambda$ is {\it completely positive} if $\Lambda \ot \Id_{\cH_E}$ is positive
for all $\cH_E$.

We seek an encoder $\cE$ and decoder $\cD$ and a compression space $\cH_{0n}$:
 $$
 \cH^{\ot n}\xrightarrow{\cE}\cH_{0n}=(\BC^2)^{\ot nR}\xrightarrow{\cD}  \cH^{\ot n}
 $$
 with $R$ as small as possible such that
 $  \cE\circ\cD(\rho^{\ot n})$ converges to $\rho^{\ot n}$ as $n\ra \infty$.
 To determine $R$, we need a quantum version of entropy.

 \begin{definition} The {\it von Neumann entropy} of  a density operator $\rho$ is
 $H(\rho)=-\ttrace(\rho\tlog(\rho))$. 
 \end{definition}
 
 Here $\tlog(\rho)$ is defined as follows:
   write $\rho$ in terms of its eigenvectors and eigenvalues, $\rho=\sum_j\l_j \ketbra{\psi_j}{\psi_j}$,
 then   $\tlog(\rho)= \sum_j \tlog(\l_j)\ketbra{\psi_j}{\psi_j}$.
 
 If $\rho=\sum_j\l_j \ketbra{\psi_j}{\psi_j}$, then
 $H(\rho)=-\sum_j \l_j\tlog(\l_j)$ so   if $\rho$ is classical (i.e., diagonal), one obtains
 the Shannon entropy.
 
 \begin{proposition}\label{vnprops} The von Neumann entropy has the following properties:
 \begin{enumerate}
 \item $H(\rho)\geq 0$ with equality if and only if $\rho$ is pure.
 \item Let $\tdim \cH=d$. Then $H(\rho)\leq \tlog(d)$ with equality if and only if $\rho=\frac 1d\Id_{\cH}$.
 \item If $\rho=\ketbra\psi\psi \in \tend(\cH_A\ot\cH_B)$, then
 $H(\rho_A)=H(\rho_B)$, where $\rho_A=\ttrace_{\cH_B}(\rho)\in \tend(\cH_A)$.
 \end{enumerate}
 \end{proposition}

\begin{theorem}\cite{MR1328824}[The quantum noiseless channel theorem]
Let $(\cH,\rho)$ be an i.i.d. quantum source. If $R>H(\rho)$,  then there exists
a reliable compression scheme of rate $R$. That is, there exists a compression space
$\cH_{0n}$, of dimension $2^{nR}$, and encoder $\cE:\cH^{\ot n}\ra \cH_{0n}$ and a decoder $\cD:\cH_{0n}\ra
\cH^{\ot n}$ such that $  \cD\circ \cE(\rho^{\ot n})$ converges to
$\rho^{\ot n}$ as $n\ra \infty$. If $R<H(\rho)$, then any compression scheme
is unreliable.
\end{theorem}
 
 \section{More on  von Neumann entropy and its variants}\label{vnplus}

 First for the classical case, define the relative entropy
 $H(\ol p|| \ol q):=-\sum p_i\tlog\frac{q_i}{p_i}=-H(\ol p)-\sum_i p_i\tlog(q_i)$.
 It is zero when $\ol p=\ol q$ and is otherwise positive.
 Define the relative von Neumann entropy $H(\rho|| \sigma):=\ttrace(\rho\tlog(\rho))-\ttrace(\rho\tlog(\s))$. 
 It shares the positivity property of its classical cousin:
  \cite{Klein} $H(\rho||\s)\geq 0$ with equality if and only if $\rho=\s$.

  von Neumann Entropy is non-decreasing under projective measurements: Let
 $\Pi_i$ be a complete set of orthogonal projectors, set $\rho'=\sum_i \Pi_i \rho\Pi_i$. Then
 $H(\rho')\geq H(\rho)$ with equality if and only if $\rho'=\rho$.
 
 \medskip

 Here and in what follows $\rho_{AB}$ is a density operator on $\cH_A\ot \cH_B$ and
 $\rho_A=\ttrace_{\cH_B}(\rho_{AB})$, $\rho_B=\ttrace_{\cH_A}(\rho_{AB})$
 are respectively the induced density operators on $\cH_A$, $\cH_B$.
 
    von Neumann entropy
    is sub-additive: 
 $H(\rho_{AB})\leq H(\rho_A)+H(\rho_B)$ with equality if and only if $\rho_{AB}=\rho_A\ot \rho_B$.
  It also satisfies a triangle inequality: $H(\rho_{AB})\geq |H(\rho_A)-H(\rho_B)|$.

\medskip

 Recall the conditional Shannon entropy is defined to be
 $H(\ol p_{\cX}|\ol p_{\cY})=-\sum_{i,j}p_{\cX\times \cY}(i,j)\tlog p_{\cX|\cY}(i|j)$,
 the entropy of $p_{\cX}$ conditioned on $y=j$, averaged over $\cY$. It is not clear how
 to \lq\lq condition\rq\rq\  one density matrix on another,  so one needs   a different
 definition. Recall that Shannon entropy satisfies
 $H(\ol p_{\cX}|\ol p_{\cY})=H(\ol p_{\cX\times \cY})-H(\ol p_{\cY})$, and the right
 hand side of this expression does make sense for density operators, so  
 define, for $\rho_{AB}$ a density operator on $\cH_A\ot \cH_B$,
 \be\label{condvn}
 H(\rho_A|\rho_B):=H(\rho_{AB})-H(\rho_B).
 \ene
 Note that $H(\rho_A|\rho_B)$ is a function of $\rho_{AB}$, as $\rho_B=\ttrace_{\cH_A}\rho_{AB}$.
 
 WARNING: it is possible that the conditional von Neumann entropy is {\it negative}
 as it is possible that $H(\rho_B)>H(\rho_{AB})$. Consider the following example:
 Let $\bra\psi=\eprx\in \cH_A\ot \cH_B$. Then $\rho_A=\frac 12\Id_{\cH_A}=
 \frac 12(\ketbra 00 +\ketbra 11)$ so $H(\rho_A)=1$, but
  $H(\ketbra\psi\psi)=0$ because $\ketbra\psi\psi$ is pure.
 
 However, vestiges of positivity are true in the quantum case:
 \begin{theorem} [Strong sub-additivity]\label{ssubadthm}
 Let $\rho_{ABC}$ be a density operator on $\cH_A\ot \cH_B\ot \cH_C$. Then
 \be\label{ss1}
 H(\rho_C|\rho_A)+H(\rho_C|\rho_B)\geq 0
 \ene
and
 \be\label{ss2}
 H(\rho_{ABC})-[H(\rho_{AB})+H(\rho_{BC})] + H(\rho_B)\geq 0.
 \ene
 \end{theorem}

 \section{Entanglement and LOCC}\label{ELOCC}
 We have seen several ways that {\it entanglement} is  a resource already for
 the space $\cH_A\ot \cH_B=\BC^2\ot \BC^2$: given a shared    $\epr=\eprx$,
 one can transport two bits of classical information using only one qubit 
 (\lq\lq super dense coding\rq\rq) and one can also transmit one qubit of 
 quantum information from Alice to Bob by sending two classical
 bits (\lq\lq teleportation\rq\rq ).

 \subsection{LOCC}Assume  several different laboratories can communicate
 classically, have prepared some shared states in advance, and can perform unitary
 and projection operations on their parts of the states, as was the situation for 
 quantum teleportation.
More precisely,    make the following assumptions:
 
 \begin{itemize}
 \item $\cH=\cH_1\otc \cH_n$, and the $\cH_j$ share an entangled
   state $\bra\psi$. Often one will just have
 $\cH=\cH_A\ot \cH_B$ and $\bra\psi=\a \bra{00} +\b\bra{11}$.
 
 \item The laboratories can communicate classically.
 
 \item Each laboratory is allowed to perform unitary and measurement operations
 on their own spaces.
 \end{itemize}
 
 The above assumptions are called {\it LOCC} for \lq\lq local operations and classical
 communication\rq\rq . It generalizes the set-up for teleportation \S\ref{teleporationsect}.
 
 Restrict to the case $\cH=\cH_A\ot \cH_B$, each of dimension two.
 I will use $\epr$  as a benchmark for measuring the quality of entanglement.
 
We will not be concerned with a single state $\bra\psi$, but the tensor
 product of many copies of it, ${\bra\psi}^{\ot n}\in (\cH_A\ot \cH_B)^{\ot n}$.
   \lq\lq How much\rq\rq\  entanglement does ${\bra\psi}^{\ot n}$ have?
An answer is given in \S\ref{ecd}.
 
 To gain insight as to  which  states can be produced via  LOCC  from a given density operator,   return to the classical case.
 For the classical cousin of LOCC, by considering diagonal density operators,
 we see we should allow alteration of a probability distribution by
 permuting the $p_j$ (permutation matrices are unitary), and more generally
 averaging
 our probability measure under some probability measure on elements of
 $\FS_d$ (the classical cousin of a projective measurement), i.e., we should allow
\be\label{classicaldegen}
 \ol{p}\mapsto
 \sum_{\s\in \FS_d}q_{\s}\mu(\s)\ol p
\ene
 where $\mu: \FS_d\ra GL_d$ is the representation, and $q$ is a probability
 distribution on $\FS_d$. 
 
 This is because the unitary and projection local operators allowed amount to
 $$
 \rho\mapsto \sum_{j=1}^k p_j U_j\rho U_j\inv
 $$
 where the $U_j$ are unitary and $p$ is a probability distribution on 
 $\{1\hd k\}$ for some finite $k$.
 
\subsection{A partial order on  probability distributions compatible with entropy}  Shannon entropy is
 non-increasing under an action of the form \eqref{classicaldegen}.
The   partial order on probability distributions determined by  \eqref{classicaldegen} 
is  
  the {\it dominance order}:
   
 \begin{definition}  Let $x,y\in \BR^d$, write $x^{\downarrow}$ for $x$ re-ordered
 such that $x_1\geq x_2\geq \cdots \geq x_d$. Write $x\prec y$ if
 for all $k\leq d$,  $\sum_{j=1}^k x_j^{\downarrow} \leq \sum_{j=1}^k y_j^{\downarrow}$.
 \end{definition}
 
 Note that if $p$ is a probability distribution concentrated at a point, then $\ol q\prec \ol p$
 for all probability distributions $q$, and if $  p$ is such that $p_j=\frac 1d$ for all $j$, then
 $\ol p\prec \ol q$ for all $q$, and more generally the dominance order is compatible with the
 entropy in the sense that $\ol p\prec \ol q$ implies $H(\ol p)\geq H(\ol q)$.

Recall that a  matrix $D\in Mat_{d\times d}$ is  doubly stochastic  if
 $D_{ij}\geq 0$ and all column and row sums equal one. Let
 $\cD\cS_d\subset Mat_{d\times d}$ denote the set of doubly stochastic matrices.
  G. Birkoff   \cite{MR0020547}  showed 
 $\cD\cS_d=\tconv(\mu(\FS_d))$, and Hardy-Littlewood-Polya \cite{MR0046395}  showed
 $\{ x\mid x\prec y\}=  \cD\cS_d\cdot y$.

 \subsection{A reduction theorem}
 The study of LOCC is potentially unwieldy because  there can be numerous
 rounds of local operations and classical communication, making it hard to model.
 The following result eliminates this problem:
 
 \begin{proposition}\label{prop5511}  If $\bra\psi\in \cH_A\ot \cH_B$ can be transformed
 into $\bra\phi$ by LOCC, then it can be transformed to $\bra\phi$ by the following
 sequence of operations:
 \begin{enumerate}
 \item Alice performs a single measurement with operators $\Pi_{M_j}$.
 \item She sends the result of her measurement (some $j$) to Bob classically.
 \item Bob performs a unitary operation on his system.
 \end{enumerate}
 \end{proposition}
  
 The key point is that for any vector spaces $V,W$,  an element $f\in V\ot W$, may be considered as
 a linear map $W^*\ra V$. In our case, $\cH_B^*\simeq \cH_B$ so $\bra\psi$ induces
 a linear map $\cH_B\ra \cH_A$ which gives us the mechanism to transfer Bob's measurements
 to Alice.

 Now I can state the main theorem on LOCC:
 
 \begin{theorem} \cite{nie99}  $\bra\psi\leadsto \bra\phi$ by LOCC
 if and only if 
 $\tsingvals(\bra\psi)\prec\tsingvals(\bra\phi)$.
 \end{theorem}

\subsection{Entanglement distillation (concentration) and dilution}\label{ecd}

To compare the entanglement resources of two states $\bra\phi$ and $\bra\psi$,
   consider ${\bra\phi}^{\ot m}$ for large $m$ with the goal of  determining the
largest $n=n(m)$ such that ${\bra\phi}^{\ot m}$ may be degenerated to ${\bra\psi}^{\ot n}$
via LOCC. Due to the approximate and probabilistic nature of quantum computing,
relax this  to degenerating ${\bra\phi}^{\ot m}$ to a state that is close to ${\bra\psi}^{\ot n}$.

There is a subtlety for this question worth pointing out. Teleportation was
defined in such a way that Alice did not need to know the state she was teleporting, but
for distillation and dilution, 
she will  need to know that its right singular vectors are standard basis vectors. More precisely,
if she is in possession of $\bra\psi= \sqrt{p_1}\bra {v_1}\ot \bra 1+
\sqrt{p_2}\bra {v_2}\ot \bra 2$ , 
she can teleport the second half of it to Bob if they share $\epr\in  \cH_A\ot \cH_B$. More generally,  
if she is in possession of 
$\bra\psi=\sum_{j=1}^d \sqrt{p_j} \bra{v_j}\ot \bra j\in \cH_{A'}\ot \cH_{A''}$, she can teleport it to 
Bob if they share enough EPR states. In most textbooks, Alice is assumed to possess states
whose singular vectors are   $\bra{jj}$'s and I will follow that convention here.
Similarly, if
$\bra\psi=\sum_{j=1}^d \sqrt{p_j} \bra{jj} \in \cH_{A }\ot \cH_{B}$,
I   discuss
how many shared EPR states they can construct from a shared $\bra\psi^{\ot m}$.

Define the {\it entanglement cost} $E_C(\psi)$ to be 
$\tinf_{m} \frac{n(m)}m$ where $n(m)$ copies of $\psi$ can be constructed from $\epr^{\ot m}$ by LOCC with
error  going to zero  as $m\ra \infty$.
Similarly, define the {\it entanglement value}, or {\it distillable entanglement}  $E_V(\psi)$ to be
$\tsup_{m} \frac{n(m)}m$ where $n(m)$ copies of $\epr$ can be constructed with diminishing error from $\bra\psi^{\ot m}$ by LOCC.
One has 
$ E_V(\psi)= E_C(\psi)= H(\ketbra\psi\psi)$.

\begin{remark} In classical computation one can reproduce information, but this cannot be done
with  quantum information in general. This is because the map $\bra\psi \mapsto \bra\psi\ot \bra\psi$,
called the {\it Veronese map} in algebraic geometry, is not a linear map.
This observation is called the {\it no cloning theorem} in the quantum literature. 
However, one can define a linear map, e.g., $\BC^2\ra \BC^2\ot \BC^2$ that duplicates basis vectors,
i.e., $\bra 0\mapsto\bra 0\ot \bra 0$ and $\bra 1\mapsto \bra 1\ot \bra 1$. But then
of course $\a\bra 0 + \b \bra 1 \mapsto \a \bra 0\ot \bra 0 + \b \bra 1\ot \bra 1\neq
(a\bra 0+\b\bra 1)^{\ot 2}$.
\end{remark}

For mixed states $\rho$ on $\cH_A\ot \cH_B$, one can still define $E_C(\rho)$ and $E_V(\rho)$, but
there exist examples where they differ, so there is not a   canonical  measure of entanglement.
A wish list of what one  might want from an entanglement measure $E$: 
\begin{itemize}
\item Non-increasing under LOCC. 
\item If $\rho$ is a product state, i.e.,  $\rho=\ketbra{\phi_A}{\phi_A}\ot \ketbra{\psi_B}{\psi_B}$, then $E(\rho)=0$.
\end{itemize}

The two conditions together imply any state constructible from 
a product state  by LOCC should also have zero
entanglement. Hence the following definition:

 \begin{definition}\label{entangledef}A density operator $\rho\in \tend (\cH_1\otc \cH_n)$ is {\it separable}
 if
 $\rho=\sum_i p_i \rho_{i,1}\otc \rho_{i,n}$, where
 $\rho_{i,\a}\in \tend(\cH_{\a})$ are density operators, $p_i\geq 0$, and
 $\sum_i p_i=n$. If $\rho$ is not separable,   $\rho$ is {\it entangled}.
 \end{definition}

 \begin{definition} An {\it entanglement monotone} $E$ is a function on density operators on $\cH_A\ot \cH_B$ that is
 non-increasing under LOCC.
 \end{definition}
 
 An example of an entanglement monotone different from $E_V,E_C$ useful for general density operators is defined in 
 \cite{MR2036165}.
 
 \section{SLOCC}\label{slocc}
An entanglement measure appealing to geometers is SLOCC (stochastic local operations and classical communication) defined originally in \cite{PhysRevA.63.012307}, which
asks if $\bra\psi\in \cH_1\otc\cH_d $ is in the same $SL(\cH_1)\ctimes SL(\cH_d)$ orbit as  $\bra\phi \in \cH_1\otc\cH_d $. 
If one relaxes this to orbit closure, then it amounts
to being able to convert $\bra\psi$ to $\bra\phi$ with positive probability. While appealing, and while there is
literature on SLOCC, given the probabilistic nature of quantum computing, its use appears to be limited
to very special cases, where the  orbit structure is understood (e.g., $d\leq 4$, $\tdim \cH_j=2$).

\section{Tensor network states}\label{tns}

Physically, entanglement is more likely when the particles are
closer together, so if we have an arrangement of electrons, say
on a circle, as in Figure \ref{figx}:

\begin{figure}[!htb]\begin{center}
\includegraphics[scale=.3]{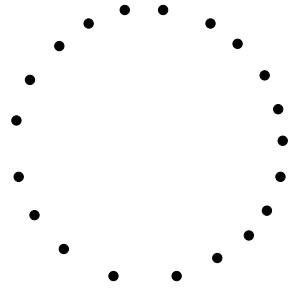}
\end{center}
\caption{\label{figx}  Electrons arranged on a circle}
\end{figure}
It is highly improbable that the electrons will share entanglement with any but their
nearest neighbors.  This is fortuitous, because if one is dealing with thousands
of electrons and would like to describe their joint state, {\it a priori} one would
have to work with a vector space of dimension $2^{n}$, with $n$ in the thousands,
 which is not feasible.
The practical solution to this problem is to define a subset of $(\BC^2)^{\ot n}$
of reasonable dimension (e.g. $O(n)$) consisting of the probable states.

\begin{figure}[!htb]\begin{center}
\includegraphics[scale=.3]{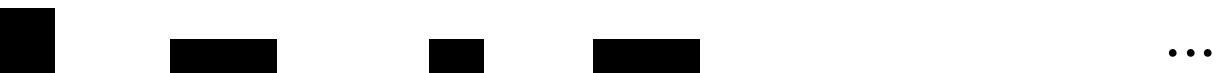}\label{figb}
\end{center}
\end{figure}

 For example, say the isolated system consists of electrons arranged along a line as in \ref{figb}.

and we only want to allow electrons to be entangled with their nearest neighbors.
This leads to the notion of {\it Matrix Product States (MPS)}: draw a graph reflecting this geometry,
with a vertex for each electron. To each vertex, attach edges going from the electron's vertex
to those of its nearest neighbors, and add an additional edge not attached to anything else
(these will be called physical edges). If our space is
$\cH_1\otc \cH_n$, then, assuming vertex $j$ has two neighbors, attach two auxiliary
vector spaces, $E_{j-1},E_{j}^*$, and a tensor $T_j\in \cH_j\ot E_{j-1}\ot E_{j}^*$.
If we are on a line, to vertex one, we just attach $T_1\in \cH_1\ot E_1^*$, and
similarly, to vertex $n$ we attach $T_n\in \cH_n\ot E_{n-1}$. Now consider the tensor
product of all the tensors
$$
T_1\otc T_n\in 
(\cH_1\ot E_1^*)\ot (\cH_2\ot E_1\ot E_2)\otc
(\cH_{n-1}\ot E_{n-2}\ot E_{n-1}^*)\ot (\cH_n\ot E_{n-1})
$$
Assume each $E_j$ has dimension $k$.
We can contract these to obtain a tensor 
$T\in \cH_1\otc \cH_n$. If $k=1$, we just obtain the product states.
As we increase $k$, we obtain a steadily larger subset of $\cH_1\otc\cH_n$.
The claim is that the tensors obtainable in this fashion (for some $k$ 
determined by    the
physical setup) are exactly those locally entangled states that we seek.

\begin{figure}[!htb]\begin{center}
\includegraphics[scale=.3]{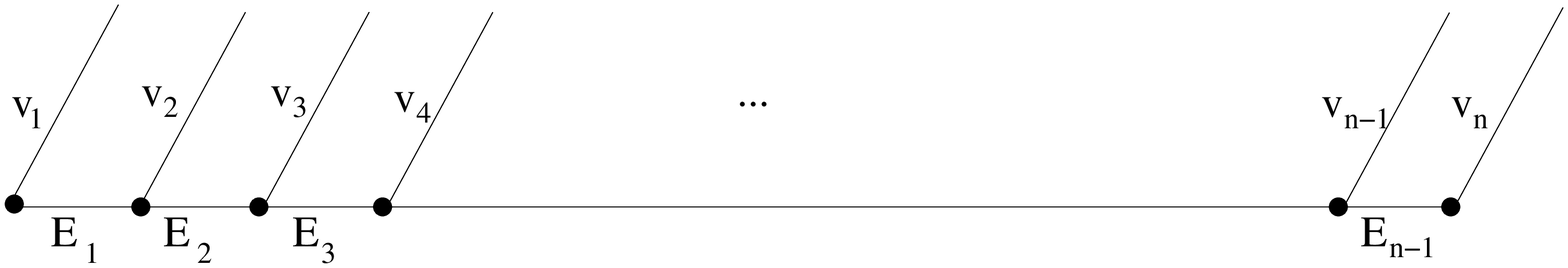}
\end{center}
\end{figure}

For the circle, the only difference in the construction  is to make the result periodic, so 
$T_1\in \cH_1\ot E_n\ot E_1^*$ and $T_n\in \cH_n\ot E_{n-1}\ot E_n^*$.

\begin{figure}[!htb]\begin{center}
\includegraphics[scale=.3]{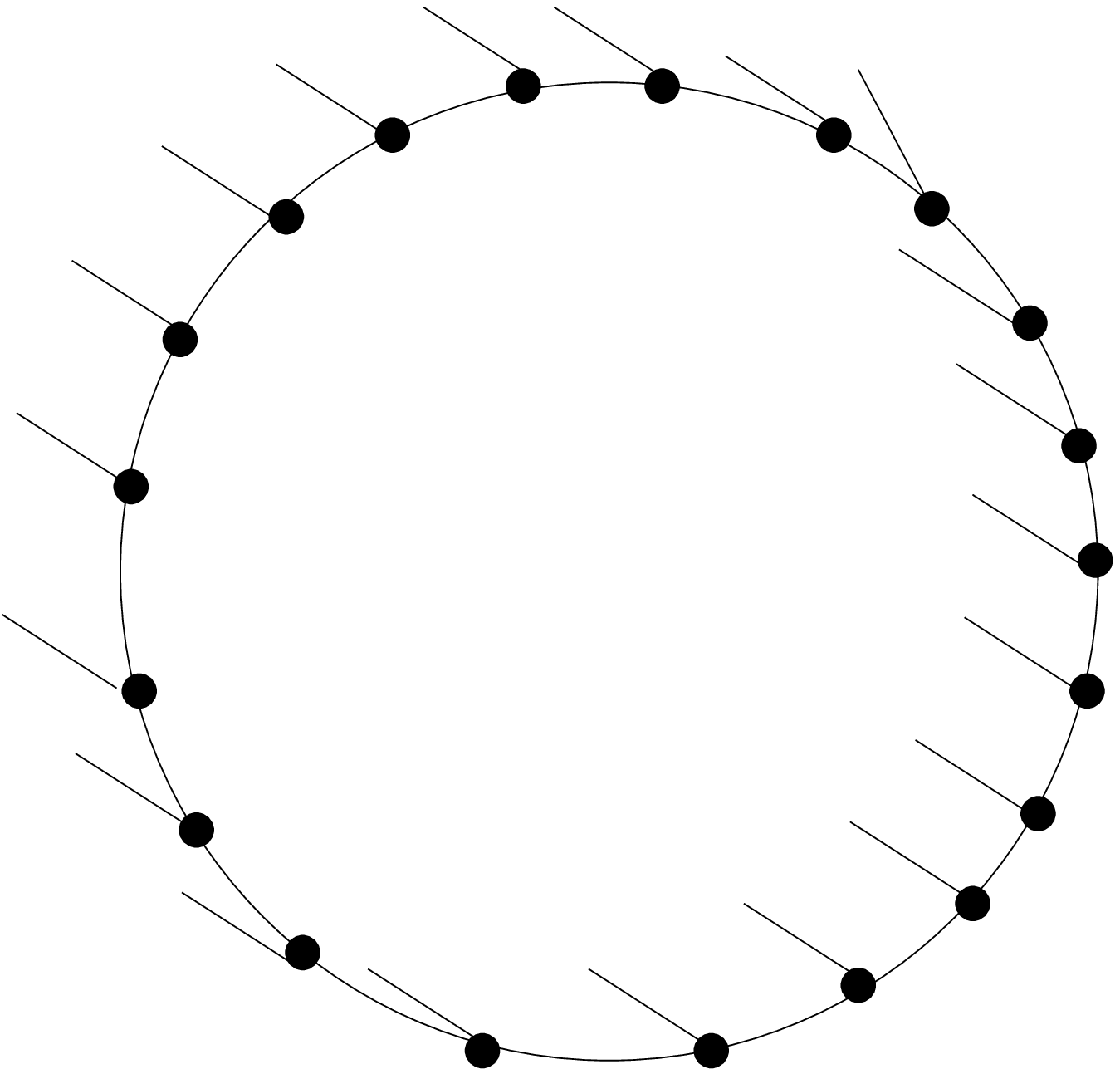}
\end{center}
\end{figure}

Sometimes for applications (e.g. translation invariant systems on the circle) one requires
the same tensor be placed at each vertex. If the tensor is $\sum_{i,j,\a}T_{i,j,\a}\ket i \ot \bra j\ot v_{\a}$,
the resulting tensor is
$\sum T_{i_1,i_2,\a_1}T_{i_2,i_3,\a_2}\cdots T_{i_n,i_1,\a_n}
v_{\a_1}\otc v_{\a_n}$. 

For  a second example,   consider electrons arranged
in a rectangular array (or on a grid on a torus), where each vertex is allowed to interact
with its four nearest neighbors.

Assume we place the same tensor at each
vertex. If our grid is $n\times n$ and periodic, 
we obtain a map $(\BC^k)^{\ot 4}\ot \BC^d\ra (\BC^d)^{\ot n^2}$.
 
\begin{definition} Let $\G$ be a directed graph with vertices $v_{\a}$ and  two kinds of edges: \lq\lq physical\rq\rq\ edges $e_i$, that
are attached to a single vertex, and \lq\lq auxiliary\rq\rq\  (or {\it entanglement}) edges $e_s$ between two vertices.
Associate to each physical edge a vector space $V_i$ (or perhaps better $\cH_i$), and to each auxiliary edge a vector space $E_s$,
of dimension $\bee_s$. Let $\ol \bee=(\bee_1\hd \bee_f)$ denote the vector of these dimensions.
A {\it tensor network state} associated to $(\G, \{ V_i\}, \ol \bee)$ is a tensor 
$T\in V_1\otc V_n$ obtained as follows: To each vertex $v_{\a}$, associate a tensor 
$$T_{\a}\in \ot_{i\in \a}V_i\ot_{s\in in(\a)}E_s^*
\ot_{t\in out(\a)}E_t.
$$
Here $in(\a)$ are the edges going into vertex $\a$ and $out(\a)$ are the edges going out of the vertex.
The {\it  tensor network state} associated to this configuration is   $T:=contr(T_1\otc T_g)\in V_1\otc V_n$.
Let $TNS(\G, V_1\otc V_n, \bee)\subset V_1\otc V_n$ denote the set of tensor network states.
\end{definition}

\begin{example} Let $\G$ be:
\begin{figure}[!htb]\begin{center}
\includegraphics[scale=.3]{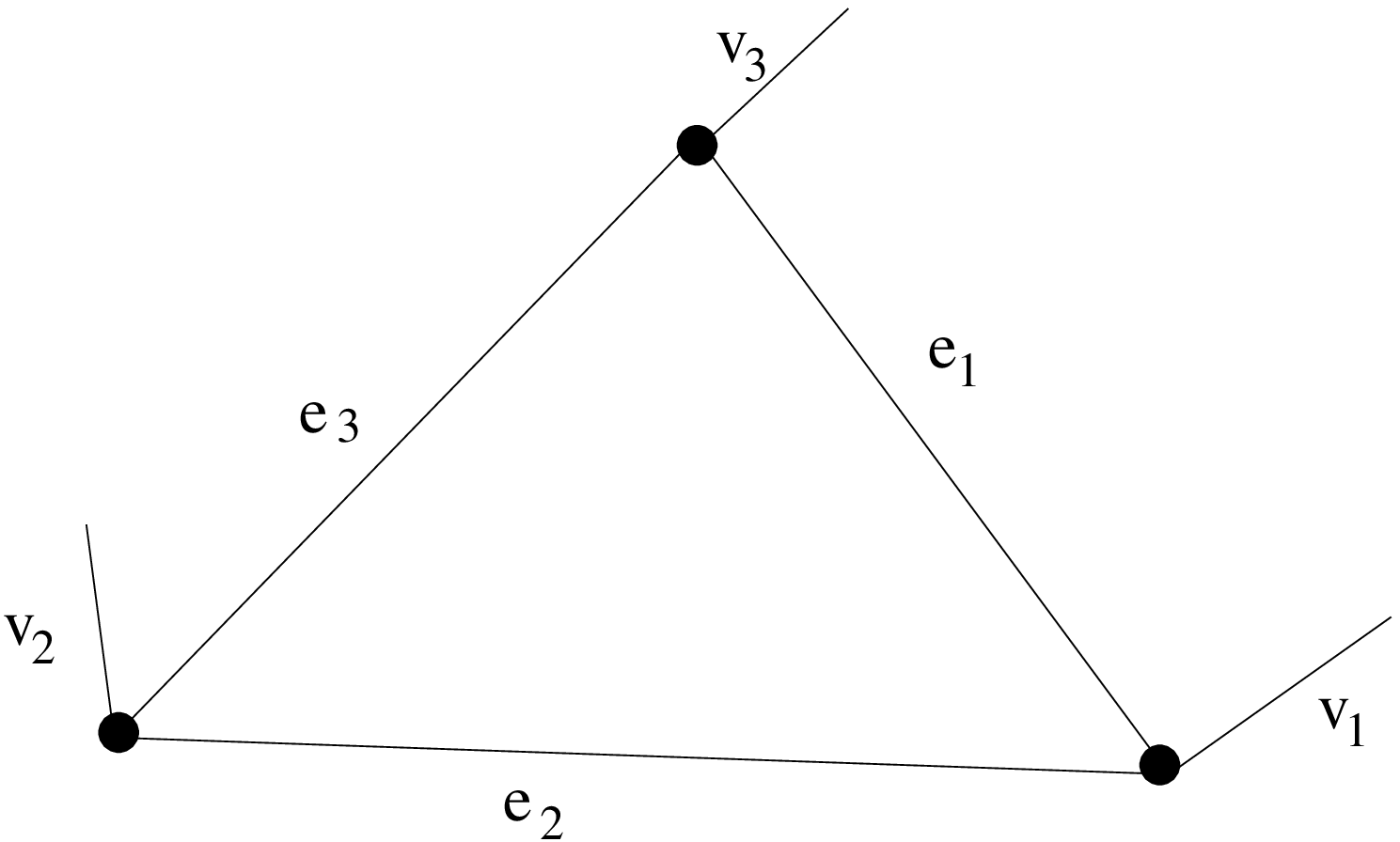} 
\end{center}
\end{figure}

Then 
\begin{align*}
TNS(\G,V_1\ot V_2\ot V_3,\ol\bee) &= TNS(\G,(E_1^*\ot E_2)\ot (E_2^*\ot E_3)\ot (E_3^*\ot E_1),\ol\bee)\\
&=
\tend(V_1)\times \tend(V_2)\times \tend(V_3)\cdot M_{\langle \bee_1,\bee_2,\bee_3\rangle}.
\end{align*}
Here $M_{\langle \bee_1,\bee_2,\bee_3\rangle}$ is the {\it matrix multiplication tensor},
for $A\in Mat_{\bee_1\times \bee_2}$, $B\in Mat_{\bee_2\times \bee_3}$, $C\in Mat_{\bee 3\times \bee 1}$, 
$(A,B,C)\mapsto \ttrace(ABC)$.
Let $e_1\hd e_{\bee_1}$ be a basis of $E_1$, $f_1\hd f_{\bee_2}$ be a basis of $E_2$,
and $g_1\hd e_{\bee_3}$ be a basis of $E_3$.
\end{example}
 
 There are many open questions about tensor network states: what are the algebraic varieties one obtains
 from the Zariski closure of a given TNS?

\section{Representation theory in quantum information theory} \label{rtinqt}

Natural projection operators on $V^{\ot d}$ are determined by representation theory. I review the relevant
representation theory and then apply it to describe the solution to the quantum marginal problem.

\subsection{Review of relevant representation theory}  
(Isomorphism classes of) irreducible representations of the permutation group $\FS_d$ are indexed
by partitions of $d$, write $[\pi]$ for the $\FS_d$-module corresponding to the partition $\pi$.
The irreducible polynomial representations of $GL(V)$ are indexed by partitions
$\pi=(p_1\hd p_{\ell(\pi)})$ with $\ell(\pi)\leq \tdim V$. Write $S_{\pi}V$ for the corresponding
$GL(V)$-module.
  
\begin{theorem} [Schur-Weyl duality] As a $GL(V)\times \FS_d$-module,
$$
V^{\ot d}=\bigoplus_{|\pi|=d} S_{\pi  }V\ot [\pi].
$$
\end{theorem}

Let  $P_{\pi}: V^{\ot d}\ra S_{\pi}V\ot[\pi]$ denote 
the $GL(V)\times \FS_d$-module projection operator.

One is often interested in decompositions of a module under the action of a
subgroup. For example $S^d(V\ot W)$ is an irreducible $GL(V\ot W)$-module, but as a
$GL(V)\times GL(W)$-module it has
the decomposition, called the {\it Cauchy formula}, 
\be\label{cauchyform}
S^d(V\ot W)=\oplus_{|\pi|=d}S_{\pi}V\ot S_{\pi}W.
\ene

We will be particularly interested in the decomposition
of $S^d(U\ot V\ot W)$ as a $GL(U)\times GL(V)\times GL(W)$-module.
An  explicit formula for this decomposition is {\it not known}.
Write
$$
S^d(U\ot V\ot W)=\bigoplus_{|\pi|,|\mu|,|\nu|=d}(S_{\pi}U\ot S_{\mu}V\ot S_{\nu}W)^{\op k_{\pi,\mu,\nu}}.
$$
The numbers $k_{\pi,\nu,\mu}$ that record the multiplicities are called {\it Kronecker coefficients}.
They have several additional descriptions.  For example,
$
S_{\pi}(V\ot W)=\bigoplus_{|\mu|,|\nu|=d} (S_{\mu}V\ot S_{\nu}W)^{\op k_{\pi,\mu,\nu}},
$
and  $k_{\pi,\mu,\nu}=\tdim ([\pi]\ot [\mu]\ot [\mu])^{\FS_d}  =\tmult([d], [\pi]\ot[\mu]\ot[\nu])=\tmult([\pi],[\mu]\ot [\nu])$.

\subsection{Quantum marginals and projections  onto isotypic subspaces of $\cH^{\ot d}$}
In this section I address the question: what are compatibility conditions on density operators
$\rho$ on $\cH_A\ot \cH_B$, $\rho'$ on $\cH_A$ and $\rho''$ on $\cH_B$ such that
$\rho'=\ttrace_{\cH_B}(\rho)$, $\rho''=\ttrace_{\cH_A}(\rho)$? As you might expect by now,
compatibility will depend only on the spectra of the operators.

Above I discussed representations of the general linear group $GL(V)$ where $V$ is a complex vector space. In quantum theory, 
one is  interested in representations on the unitary group $\bU(\cH)$ on a Hilbert space $\cH$. 
The unitary group is a real Lie group, not a complex Lie group, because
complex conjugation is not a complex linear map.
It is a special case of a general fact about representations of a maximal   compact subgroups of complex Lie groups
have the same representation theory as the the original group, so in particular
the decomposition of $\cH^{\ot d}$ as a $\bU(\cH)$-module coincides with its decomposition as a $GL(\cH)$-module.

For a partition $\pi=(p_1\hd p_d)$ of $d$, introduce the notation $\ol \pi=(\frac {p_1}d\hd \frac{p_d}d)$ which is 
a probability distribution on $\{1\hd d\}$.

\begin{theorem}\cite{MR2197548} \label{cmmarginal} Let $\rho_{AB}$ be a density operator
on $\cH_A\ot \cH_B$. Then there exists a sequence
$(\pi_j,\mu_j,\nu_j)$ of triples of partitions such that
$k_{\pi_j,\mu_j,\nu_j}\neq 0$ for all $j$ and 
\begin{align*}
\tlim_{j\ra \infty}\ol{\pi}_j &=\tspec(\rho_{AB})\\
\tlim_{j\ra \infty}\ol{\mu}_j &=\tspec(\rho_{A})\\
\tlim_{j\ra \infty}\ol{\nu}_j &=\tspec(\rho_{B}).\end{align*}
\end{theorem}

\begin{theorem}\cite{Klyachkopreprint} \label{Kmarginal1}
Let $\rho_{AB}$ be a density operator
on $\cH_A\ot \cH_B$ such that $\tspec(\rho_{AB})$, $\tspec(\rho_A)$ and $\tspec(\rho_B)$ are all
rational vectors. Then there exists an integer  $M>0$ such that
$$k_{M\tspec(\rho_A),M\tspec(\rho_B),M\tspec(\rho_C)}\neq 0.
$$
\end{theorem}

\begin{theorem}\cite{Klyachkopreprint} \label{Kmarginal2} Let $\pi,\mu,\nu$ be partitions of $d$
with $k_{\pi,\mu,\nu}\neq 0$ and
satisfying $\ell(\pi)\leq mn$, $\ell(\mu)\leq m$, and $\ell(\nu)\leq n$. Then there exists
a density operator $\rho_{AB}$ on $\BC^n\ot \BC^m=\cH_A\ot \cH_B$ with
$\tspec(\rho_{AB})=\ol\pi$, $\tspec(\rho_{A })=\ol\mu$, and $\tspec(\rho_{B})=\ol\nu$.
\end{theorem}

Klyatchko's proofs are via co-adjoint orbits and vector bundles on flag varieties, while the proof  of Christandl-Mitchison is information-theoretic in 
flavor.  

Recall the relative entropy $H(\ol p||\ol q)=-\sum_i p_i\tlog\frac{q_i}{p_i}$, which may be thought of as measuring
how close $p,q$ are because it is non-negative, and zero if and only if $p=q$.
A key step in the Christandl-Mitchison proof is the following theorem: 
\begin{theorem}\label{KWthm}\cite{MR1878924}
Let $\rho\in \tend(\cH)$ be a density operator, where $\tdim \cH=n$.
Let $|\pi|=d$ and let  $P_{\pi}: \cH^{\ot d}\ra S_{\pi}\cH\ot [\pi]$ be the projection operator. Then
$$
\ttrace(P_{\pi}\rho^{\ot d})\leq (d+1)^{\binom n2}e^{-d H(\ol\pi ||\tspec(\rho))}.
$$
\end{theorem}
 A key step of the proof is that   the projection of $e_I$ to $S_{\pi}V\ot [\pi]$ is nonzero
if and only if $wt(e_I)\prec \pi$.

Note that one can use quantum theory to deduce representation-theoretic consequences:
$k_{\mu,\nu,\pi}\neq 0$ implies
$H(\ol \pi)\leq H(\ol\mu)+H(\ol\nu)$, $H(\ol \mu)\leq H(\ol\pi)+H(\ol\nu)$,
and $H(\ol \nu)\leq H(\ol\mu)+H(\ol\pi)$.

Let  $Spec_{m,n,mn}$ denote the set of admissible triples $(\tspec(\rho_A),\tspec(\rho_B),\tspec(\rho_{AB})$ 
and $KRON_{m,n,mn}$ the triples $(\ol\mu,\ol\nu,\ol\pi)$ of normalized partitions $(\mu,\nu,\pi)$ with $\ell(\mu)\leq m$, $\ell(\nu)\leq n$, $\ell(\pi)\leq mn$
and $k_{\pi,\mu,\nu}\neq 0$.

The theorems above imply:
$$
Spec_{m,n,mn}=\ol{KRON_{m,n,mn}}.
$$
In particular, $Spec_{m,n,mn}$ is a convex polytope.

\subsection*{Acknowledgements} I thank the organizers of the International workshop on 
Quantum Physics and Geometry, especially Alessandra Bernardi, who also co-organized an intensive summer class on
Quantum computing and quantum information theory that I gave June-July 2017. I also thank L. Chiantini, F. Gesmundo, F. Holweck, and G. Ottaviani for useful
comments on a draft of this article.

\bibliographystyle{amsalpha}
\bibliography{Lmatrix}
\end{document}